\documentclass{amsart}

\usepackage{amsmath} 
\usepackage{amssymb}

\newtheorem{theorem}{Theorem}[section] 
\newtheorem{claim}{Claim}[theorem]
\newtheorem{lemma}[theorem]{Lemma} 
\newtheorem{proposition}[theorem]{Proposition} 
\newtheorem{corollary}[theorem]{Corollary} 

\theoremstyle{definition}
\newtheorem{definition}[theorem]{Definition}

\theoremstyle{remark}
\newtheorem{remark}[theorem]{Remark}

\numberwithin{equation}{section}
\newcommand{\forces}{\Vdash} 
\newcommand{\bV}{{\bf V}} 
\newcommand{\comp}{\circ} 
\newcommand{\con}{{\mathfrak c}} 
\newcommand{\starQ}{{\bQ^\mtree_4(K^*,\Sigma^*,\bF^*)}}
\newcommand{\rQ}{{\bQ^\mtree_4(K^r,\Sigma^r,\bF^r)}}

\newcommand{\can}{2^{\textstyle \omega}} 
\newcommand{\fs}{2^{\textstyle <\!\omega}} 
\newcommand{\baire}{\omega^{\textstyle \omega}}

\newcommand{\conc}{{}^\frown\!}
\newcommand{\lh}{{\rm lh}} 
\newcommand{\rest}{{\restriction}}
\newcommand{\mrot}{{\rm root}\/} 
\newcommand{\suc}{{\rm succ}} 

\newcommand{\dom}{{\rm dom}} 
\newcommand{\rng}{{\rm rng}}

\newcommand{\cf}{{\rm cf}}
\newcommand{\nor}{{\rm {\bf nor}}\/} 
\newcommand{\pos}{{\rm pos}}
\newcommand{\val}{{\bf val}}
\newcommand{\dis}{{\bf dis}}
\newcommand{\rht}{{\rm ht}}
\newcommand{\otp}{{\rm otp}}
\newcommand{\FC}{{\rm FC}}
\newcommand{\TCR}{{\rm LTCR}}
\newcommand{\tree}{{\rm tree}}
\newcommand{\mtree}{{\rm mt}}

\newcommand{\cA}{{\mathcal A}}

\newcommand{\bc}{{\bf c}}

\newcommand{\bH}{{\bf H}}
\newcommand{\bh}{{\bf h}}
\newcommand{\bF}{{\bf F}}

\newcommand{\bk}{{\bf k}}

\newcommand{\cN}{{\mathcal N}}
\newcommand{\bP}{{\mathbb P}}
\newcommand{\cP}{{\mathcal P}}

\newcommand{\bQ}{{\mathbb Q}}

\newcommand{\mbR}{{\mathbb R}}

\newcommand{\bs}{{\bf s}}

\newcommand{\bt}{{\bf t}}

\newcommand{\cY}{{\mathcal Y}}
\newcommand{\cZ}{{\mathcal Z}}

\begin{document}

\title{Measured creatures}

\author{Andrzej Ros{\l}anowski}
\address{Department of Mathematics\\
 University of Nebraska at Omaha\\
 Omaha, NE 68182-0243, USA\\
 and Mathematical Institute of Wroclaw University\\
 50384 Wroclaw, Poland} 
\email{roslanowski@unomaha.edu}
\urladdr{http://www.unomaha.edu/$\sim$aroslano}
\thanks{The first author thanks the Hebrew University of Jerusalem for
  support during his visits to Jerusalem and the KBN
  (Polish Committee of Scientific Research) for partial support through
  grant 2P03A03114.} 

\author{Saharon Shelah}
\address{Institute of Mathematics\\
 The Hebrew University of Jerusalem\\
 91904 Jerusalem, Israel\\
 and  Department of Mathematics\\
 Rutgers University\\
 New Brunswick, NJ 08854, USA}
\email{shelah@math.huji.ac.il}
\urladdr{http://www.math.rutgers.edu/$\sim$shelah}
\thanks{The research of the second author was partially supported by the
 Israel Science Foundation. Publication 736} 

\subjclass{Primary 03E35; Secondary 03E75, 28A20, 54H05}

\begin{abstract}
We prove that two basic questions on outer measure are undecidable. 
First we show that consistently 
\begin{itemize}
\item every sup-measurable function $f:\mbR^2\longrightarrow\mbR$ is
measurable.
\end{itemize}
The interest in sup-measurable functions comes from differential
equations and the question for which functions $f:\mbR^2\longrightarrow
\mbR$ the Cauchy problem
\[y'=f(x,y),\qquad y(x_0)=y_0\]
has a unique almost-everywhere solution in the class $AC_l(\mbR)$ of
locally absolutely continuous functions on $\mbR$.

Next we prove that consistently
\begin{itemize}
\item every function $f:\mbR\longrightarrow\mbR$ is continuous on some
set of positive outer Lebesgue measure.
\end{itemize}
This says that in a strong sense the family of continuous functions (from
the reals to the reals) is dense in the space of arbitrary such functions. 

For the proofs we discover and investigate a new family of nicely definable
forcing notions (so indirectly we deal with nice ideals of subsets of the 
reals -- the two classical ones being the ideal of null sets and the
ideal of meagre ones). 

Concerning the method, i.e., the development of a family of forcing
notions (equivalently, nice ideals), the point is that whereas there
are many such objects close to the Cohen forcing (corresponding to the ideal
of meagre sets), little has been known on the existence of relatives of
the random real forcing (corresponding to the ideal of null sets),
and we look exactly at such forcing notions.
\end{abstract}

\maketitle

\section{Introduction}
The present paper deals with two, as it occurs closely related, problems
concerning real functions. The first one is the question if it is possible
that all superposition--measurable functions are measurable.
\begin{definition}
\label{supmeas}
A function $f:\mbR^2\longrightarrow\mbR$ is {\em superposition--measurable
\/} (in short: {\em sup--measurable\/}) if for every Lebesgue measurable
function $g:\mbR\longrightarrow\mbR$ the superposition 
\[f_g:\mbR\longrightarrow\mbR:x\mapsto f(x,g(x))\]
is Lebesgue measurable.
\end{definition}
The interest in sup-measurable functions comes from differential equations
and the question for which functions $f:\mbR^2\longrightarrow\mbR$ the
Cauchy problem 
\[y'=f(x,y),\qquad y(x_0)=y_0\]
has a unique almost-everywhere solution in the class $AC_l(\mbR)$ of locally
absolutely continuous functions on $\mbR$. For the detailed discussion of
this area we refer the reader to Balcerzak \cite{Ba91b}, Balcerzak and
Ciesielski \cite{BaCi97} and Kharazishvili \cite{Kh97}. Grande and
Lipi\'nski \cite{GL78} proved that, under CH, there is a non-measurable
function which is sup-measurable. Later the assumption of CH was weakened
(see Balcerzak \cite[Thm 2.1]{Ba91b}), however the question if one can build
a non-measurable sup-measurable function in ZFC remained open (it was
formulated in Balcerzak \cite[Problem 1.10]{Ba91b} and Ciesielski
\cite[Problem 5]{Ci97}, and implicitly in Kharazishvili \cite[Remark
4]{Kh97}). In the third section we will answer this question by showing
that, consistently, every sup-measurable function is Lebesgue measurable.

Next we deal with von Weizs\"acker's problem. It has enjoyed considerable
popularity, and it has origins in measure theory and topology. In
\cite{We79}, von Weizs\"acker noted that if   
\begin{enumerate}
\item[$(*)$] ${\rm non}(\cN)\stackrel{\rm def}{=}\min\{|X|:X\subseteq\mbR$
has positive outer Lebesgue measure$\;\}=\con$,
\end{enumerate}
then 
\begin{enumerate}
\item[$(\otimes)$] there is a function $f:[0,1]\longrightarrow [0,1]$ such
that the graph of $f$ is of (two dimensional) outer measure 1 but for every
Borel function $g:[0,1]\longrightarrow [0,1]$ the set $\{x\in [0,1]:f(x)=
g(x)\}$ is of measure zero.
\end{enumerate}
Then he showed that $(\otimes)$ implies
\begin{enumerate}
\item[$(\boxtimes)$] there is a countably generated
$\sigma$--algebra $\cA$ containing ${\rm Borel}([0,1])$ such that the
Lebesgue measure can be extended to $\cA$, but there is no extremal
extension to $\cA$.  
\end{enumerate}
So it was natural to ask if the statement in $(\otimes)$ can be proved in
ZFC (i.e., without assuming $(*)$). A way to formulate this question was to
ask 
\begin{itemize}
\item[$(\circledast)_{\rm vW}$] Is it consistent to suppose that for every
function $f:\mbR\longrightarrow\mbR$ there is a Borel measurable function
$g:\mbR\longrightarrow\mbR$ such that the set $\{x\in\mbR: f(x)=g(x)\}$
is not Lebesgue negligible ?
\end{itemize}

One can arrive to question $(\circledast)_{\rm vW}$ also from the
topological side. In \cite{Bl22}, Blumberg proved that if $X$ is a separable 
complete metric space and $f:X\longrightarrow\mbR$, then there exists a
dense (but possibly countable) subset $D$ of $X$ such that the restriction
$f\restriction D$ is continuous.  This result has been generalized in many
ways: by considering functions on other topological spaces, or by aiming at
getting ``a large set'' on which the function is continuous. For example, in
the second direction, we may restrict ourselves to $X=\mbR$ and ask if above
we may request that the set $D$ is uncountable. That was answered by
Abraham, Rubin and Shelah  who showed in \cite{ARSh:153} that, consistently,
every real function is continuous on an uncountable set. The next natural
step is to ask if we can demand that the set $D$ is of positive outer
measure, and this is von Weizs\"acker's question $(\circledast)_{\rm
vW}$. It appears in Fremlin's list of problems as \cite[Problem AR(a)]{Fe94}
and in Ciesielski \cite[Problem 1]{Ci97}. 

We will answer question $(\circledast)_{\rm vW}$ in affirmative in the
fourth section. The respective model is built by a small modification of the
iteration used to deal with the sup-measurability problem (and, as a matter
of fact, it may serve for both purposes). We do not know if $(\boxtimes)$
fails in our model (and the question if $\neg (\boxtimes)$ is consistent
remains open). 

Let us note that the close relation of the two problems solved here is
not very surprising. Some connections were noticed already in Balcerzak and
Ciesielski \cite{BaCi97}. Also, among others, these connections motivated
the following strengthening of $(\circledast)_{\rm vW}$:
\begin{enumerate}
\item[$(\circledast)_{\rm vW}^+$] Is it consistent that for every subset $Y$
of $\mbR$ of positive outer measure and every function $f:Y\longrightarrow
\mbR$, there exists a set $X\subseteq Y$ of positive outer measure such that
$f\restriction Y$ is continuous? 
\end{enumerate}
However, as Fremlin points out, the answer to $(\circledast)_{\rm vW}^+$ is
NO: 

\begin{proposition}
[Fremlin \cite{Fr00}]
There are a set $Y\subseteq \mbR$ of positive outer measure and a function
$f:Y\longrightarrow\mbR$ such that $f\restriction X$ is not continuous for 
any $X\subseteq Y$ of positive outer measure.
\end{proposition}

\begin{proof}
Recall that a Hausdorff space $Z$ is {\em universally negligible\/} if there
is no Borel probability measure on $Z$ that vanishes at singletons. By
Grzegorek \cite{Gr81}, there is a universally negligible set
$Z\subseteq\mbR$ of cardinality ${\rm non}(\cN)$ (see also \cite[Volume IV,
439E(c)]{Fr0x}). Pick a non-null set $Y\subseteq\mbR$ of size ${\rm
non}(\cN)$ and fix a bijection $f:Y\longrightarrow Z$. 

If $X\subseteq Y$ is such that $f\restriction X$ is continuous, then we may
transport Borel measures on $X$ to $Z$, and therefore $X$ is universally
negligible and thus Lebesgue negligible. (See also \cite[Volume IV,
439C(f)]{Fr0x}.) 
\end{proof}

The notion of sup-measurability has its category version (defined naturally 
by replacing ``Lebesgue measurability'' by ``Baire property''). It was
investigated in E.Grande and Z.Grande \cite{GG84}, Balcerzak \cite{Ba91b}, 
and Ciesielski and Shelah \cite{CiSh:695}. The latter paper presents a 
model in which every Baire-sup-measurable function has the Baire
property. Also von Weizs\"acker problem has its category counterpart which
was answered in Shelah \cite{Sh:473}. What is somewhat surprising, is that
the models of \cite{CiSh:695} and \cite{Sh:473} seem to be totally
unrelated (while for the measure case presented here the connection is
striking). Moreover, neither the forcing used in \cite{CiSh:695} (based on
the oracle-cc method of Shelah \cite[Chapter IV]{Sh:f}), nor the one applied
in Shelah \cite{Sh:473}, are parallel to the method presented here. 
\medskip

The present paper is a part of author's program to investigate the family of
forcing notions with {\em norms on possibilities}. One of the points is that
we know many forcing notions in the neighbourhood of the Cohen forcing
notion (see, e.g., Ros{\l}anowski and Shelah \cite{RoSh:628},
\cite{RoSh:672}), but we have not known any relatives of the random real
forcing. In the present paper we further develop the theory of forcing
notions with norms on possibilities introducing {\em measured
creatures}. This enrichment of the method of norms on possibilities creates 
a bridge between the forcings of \cite{RoSh:470} and the random real forcing
(including the latter in our framework), and we come with $\baire$--bounding
friends of the random forcing. Though they are not ccc, they do make random
not so lonely.  
\medskip

Our presentation is self-contained, and though we use the notation of
\cite{RoSh:470}, the two basic definitions we need from there are stated 
in somewhat restricted form below (in \ref{treecreature}, \ref{treeforc}). 
The general construction of forcing notions using measure (tree) creatures
is presented in the first section, and only in the following section we
define the particular example that works for us. The forcing notion $\starQ$
(defined in section 2) is the basic ingredient of our construction. The
required models are obtained by CS iterations of $\starQ$; in the fourth
section we also add in the iteration random reals (on a stationary set of
coordinates).   

Let us point out that ``measured creatures'' presented here have their ccc
relative which appeared in \cite[\S 2.1]{RoSh:672} (and more general
constructions will be presented in \cite[Chapter 2]{RoSh:670}).  
\medskip

\noindent{\bf Notation:}\qquad 
Most of our notation is standard and compatible with that of classical
textbooks on Set Theory (like Bartoszy\'nski and Judah \cite{BaJu95}). 
However in forcing we keep the convention that {\em a stronger condition is
the larger one}.   

\begin{enumerate}
\item $\mbR^{{\geq}0}$ stands for the set of non-negative reals. For a real
number $r$ and a set $A$, the function with domain $A$ and the constant
value $r$ will be denoted $r_A$. 
\item For two sequences $\eta,\nu$ we write $\nu\vartriangleleft\eta$
whenever $\nu$ is a proper initial segment of $\eta$, and $\nu
\trianglelefteq\eta$ when either $\nu\vartriangleleft\eta$ or $\nu=\eta$. 
The length of a sequence $\eta$ is denoted by $\lh(\eta)$.
\item A {\em tree} is a family $T$ of finite sequences such that for some
$\mrot(T)\in T$ we have
\[(\forall\nu\in T)(\mrot(T)\trianglelefteq \nu)\quad\mbox{ and }\quad
\mrot(T)\trianglelefteq\nu\trianglelefteq\eta\in T\ \Rightarrow\ \nu\in T.\]
\item For a tree $T$, the family of all $\omega$--branches through $T$ is
denoted by $[T]$, and we let 
\[\max(T)\stackrel{\rm def}{=}\{\nu\in T:\mbox{ there is no }\rho\in
T\mbox{ such that }\nu\vartriangleleft\rho\}.\]   
If $\eta$ is a node in the tree $T$ then 
\[\begin{array}{lcl}
\suc_T(\eta)&=&\{\nu\in T: \eta\vartriangleleft\nu\ \&\ \lh(\nu)=\lh(\eta)+1
\}\ \mbox{ and}\\
T^{[\eta]}&=&\{\nu\in T:\eta\trianglelefteq\nu\}.
  \end{array}\]
A set $F\subseteq T$ is {\em a  front of $T$} if 
\[(\forall \eta\in [T])(\exists k\in\omega)(\eta\restriction k\in F).\] 
\item The Cantor space $\can$ (the spaces of all functions from $\omega$ to
$2$) and the space $\prod\limits_{i<\omega} N_i$ (where $N_i$ are positive
integers thought of as non-empty finite sets) are equipped with natural
(Polish) topologies, as well with as with standard product measure
structures.  
\item For a forcing notion $\bP$, $\Gamma_\bP$ stands for the canonical
$\bP$--name for the generic filter in $\bP$. With this one exception, all
$\bP$--names for objects in the extension via $\bP$ will be denoted with a
dot above (e.g.~$\dot{\tau}$, $\dot{X}$).
\item For a relation $R$ (a set of ordered pairs), $\rng(R)$ and $\dom(R)$
stand for the range and the domain of $R$, respectively.
\item We will keep the convention that $\sup(\emptyset)$ is $0$. Similarly,
the sum over an empty set of reals is assumed to be 0.
\end{enumerate}

Let us recall the definition of tree creating pairs. Since we are going to
use local tree creating pairs only, we restrict ourselves to this case. For
more information and properties of tree creating pairs and related forcing
notions we refer the reader to \cite[\S 1.3, 2.3]{RoSh:470}.   

\begin{definition}
\label{treecreature}
Let $\bH$ be a function with domain $\omega$.
\begin{enumerate}
\item {\em A local tree--creature for $\bH$\/} is a triple 
\[t=(\nor,\val,\dis)=(\nor[t],\val[t],\dis[t])\]
such that $\nor\in\mbR^{{\geq}0}$, $\dis\in {\mathcal H}(\aleph_1)$ (i.e.,
$\dis$ is hereditarily countable), and for some sequence $\eta\in
\prod\limits_{i<n}\bH(i)$, $n<\omega$, we have 
\[\emptyset\neq\val\subseteq\{\langle \eta,\nu\rangle: \eta\vartriangleleft
\nu\in \prod_{i\leq\lh(\eta)}\bH(i)\}.\] 
For a tree--creature $t$ we let $\pos(t)\stackrel{\rm def}{=}\rng(\val[t])$.

The set of all local tree--creatures for $\bH$ will be denoted by
$\TCR[\bH]$, and for $\eta\in\bigcup\limits_{n<\omega}\prod\limits_{i<n}
\bH(i)$ we let $\TCR_\eta[\bH]=\{t\in\TCR[\bH]:\dom(\val[t])=\{\eta\}\}$. 
\item Let $K\subseteq\TCR[\bH]$. We say that a function $\Sigma:K
\longrightarrow\cP(K)$ is a {\em local tree composition on $K$\/} whenever
the following conditions are satisfied. 
\begin{enumerate}
\item[(a)] If $t\in K\cap\TCR_\eta[\bH]$, $\eta\in\prod\limits_{i<n}\bH(i)$, 
$n<\omega$, then $\Sigma(t)\subseteq\TCR_\eta[\bH]$. 
\item[(b)] If $s\in\Sigma(t)$, then $\val[s]\subseteq\val[t]$.
\item[(c)] {[{\em transitivity\/}]} If $s\in\Sigma(t)$, then $\Sigma(s)
\subseteq\Sigma(t)$.
\end{enumerate}
\item If $K\subseteq\TCR[\bH]$ and $\Sigma$ is a local tree composition
operation on $K$, then $(K,\Sigma)$ is called {\em a local tree--creating 
pair for $\bH$}. 
\item We say that $(K,\Sigma)$ is {\em strongly finitary\/} if $\bH(m)$ is
finite (for $m<\omega$) and $\TCR_\eta[\bH]\cap K$ is finite (for each
$\eta$). 
\end{enumerate}
\end{definition}

\begin{definition}
[See {\cite[Definition 1.3.5]{RoSh:470}}]
\label{treeforc}
Let $(K,\Sigma)$ be a local tree--creating pair for $\bH$. The forcing
notion $\bQ^{\tree}_4(K,\Sigma)$ is defined as follows.   
\medskip

\noindent {\bf A condition} is a system $p=\langle t_\eta:\eta\in T\rangle$ 
such that 
\begin{enumerate}
\item[(a)] $T\subseteq\bigcup\limits_{n\in\omega}\prod\limits_{i<n}\bH(i)$ is
a non-empty tree with $\max(T)=\emptyset$,
\item[(b)] $t_\eta\in\TCR_\eta[\bH]\cap K$ and $\pos(t_\eta)=\suc_T(\eta)$,
\item[(c)$_4$] for every $n<\omega$, the set 
\[\{\nu\in T:(\forall\rho\in T)(\nu\vartriangleleft\rho\ \Rightarrow\
\nor[t_\rho]\geq n)\}\]
contains a front of the tree $T$.
\end{enumerate}

\noindent{\bf The order} is given by:

\noindent $\langle t^1_\eta: \eta\in T^1\rangle\leq\langle t^2_\eta:\eta\in
T^2\rangle$ if and only if 

\noindent $T^2\subseteq T^1$ and $t^2_\eta\in\Sigma(t^1_\eta)$ for each
$\eta\in T^2$.  
\medskip

\noindent If $p=\langle t_\eta:\eta\in T\rangle$, then we write
$\mrot(p)=\mrot(T)$, $T^p= T$, $t^p_\eta = t_\eta$ etc.

The forcing notion $\bQ^\tree_\emptyset(K,\Sigma)$ is defined similarly, but
we omit the norm requirement (c)$_4$. (So $\bQ_\emptyset^\tree(K,\Sigma)$ is
trivial in a sense; we will use it for notational convenience only.)
\end{definition}

\section{Measured Creatures}
Below we introduce a relative of the {\em mixtures with random} presented in  
\cite[\S 2.1]{RoSh:672}. Here, however, the interplay between the norm of a
tree creature $t$, the set of possibilities $\pos(t)$ and the averaging
function $F_t$ assigned to $t$ is different. A more general variant of this
method will be presented in \cite[\S 2]{RoSh:670} (where we will also
consider averaging functions with discrete ranges, and also we will allow 
more complicated ``error terms'' in the equations and inequalities we will
require from them). 
\medskip

\noindent {\bf Basic Notation:} In this section, $\bH$ stands for a function 
with domain $\omega$ and such that $(\forall m\in\omega)(|\bH(m)|\geq 2)$.
Moreover we demand $\bH\in {\mathcal H}(\aleph_1)$ (i.e., $\bH$ is
hereditarily countable).

\begin{definition}
\label{mixing}
\begin{enumerate}
\item {\em A measured (tree) creature for $\bH$} is a pair $(t,F_t)$ such
that $t\in\TCR[\bH]$ and
\[F_t:[0,1]^{\textstyle \pos(t)}\longrightarrow[0,1].\]
\item We say that $(K,\Sigma,\bF)$ is a {\em measured tree creating triple
for $\bH$} if     
\begin{enumerate}
\item[(a)] $(K,\Sigma)$ is a local tree--creating pair for $\bH$, 
\item[(b)] $\bF$ is a function with domain $K$, $\bF:t\mapsto F_t$, such
that $(t,F_t)$ is a measured (tree) creature (for each $t\in K$).
\end{enumerate}
\item If $(K,\Sigma,\bF)$ is as above, $t\in K$, $X\subseteq\pos(t)$, and
$\langle r_\nu:\nu\in X\rangle\in [0,1]^{\textstyle X}$, then we define
$F_t(r_\nu:\nu\in X)$ as $F_t(r^*_\nu: \nu\in\pos(t))$, where 
\[r^*_\nu=\left\{\begin{array}{ll}
r_\nu&\mbox{ if }\nu\in X,\\
0    &\mbox{ if }\nu\in\pos(t)\setminus X.
		 \end{array}\right.\]
\end{enumerate}
\end{definition}

We think of $F_t$ as a kind of averaging function. At the first reading the 
reader may think that $\pos(t)$ is finite and 
\[F_t(r_\nu:\nu\in\pos(t))=\frac{\sum\{r_\nu:\nu\in\pos(t)\}}{|\pos(t)|}.\]
For this particular function, our construction results in the random real
forcing. However in general our averaging function does not have to be
additive (as long as it has the properties stated in \ref{nicemix} below),
and the result is not the random forcing (and this is one of the points of
our construction). Also having $F_t$ depend on $t$ allows us to ``cheat'':
if we do not like the results of our averaging we may pass to a tree
creature $s\in\Sigma(t)$ (dropping the norm a little) with possible better
for us averaging function $F_s$. 

Regarding the requirements of \ref{nicemix}, note that they are meant to
provide us with some features of the Lebesgue measure, without imposing
additivity on the averaging functions $F_t$ (specifically see
\ref{nicemix}$(\beta)$).  

\begin{definition}
\label{nicemix}
A measured tree creating triple $(K,\Sigma,\bF)$ is {\em nice} if for every
$t\in K$:
\begin{enumerate} 
\item[$(\alpha)$] if $\langle r_\nu:\nu\in\pos(t)\rangle,\langle
r_\nu':\nu\in\pos(t)\rangle\subseteq [0,1]$, $r_\nu\leq r_\nu'$ for all
$\nu\in\pos(t)$, then  
\[F_t(r_\nu:\nu\in\pos(t))\leq F_t(r_\nu':\nu\in\pos(t)),\]
\item[$(\beta)$] if $\nor[t]>1$, $\{\eta\}=\dom(\val[t])$, $r_\nu,r^0_\nu,
r^1_\nu\in [0,1]$ (for $\nu\in\pos(t)$) are such that $r^0_\nu+r^1_\nu
\geq r_\nu$ and $F_t(r_\nu:\nu\in\pos(t))\geq 2^{-2^{\lh(\eta)}}$, then
there are real numbers $c_0,c_1$ and tree creatures $s_0,s_1\in\Sigma(t)$
such that 
\[c_0+c_1=(1-2^{-2^{\lh(\eta)}})F_t(r_\nu:\nu\in\pos(t))\]
and 
\begin{enumerate}
\item[$(\otimes)$] if $\ell<2$, $c_\ell>0$, then $\nor[s_\ell]\geq\nor[t]-
1$, $\pos(s_\ell)\subseteq \{\nu\in\pos(t):r^\ell_\nu>0\}$, and 
\[F_{s_\ell}(r_\nu^\ell:\nu\in\pos(s_\ell))\geq c_\ell,\]
\end{enumerate}
\item[$(\gamma)$] if $b\in [0,1]$ and $r_\nu\in [0,1]$ (for
$\nu\in\pos(t)$), then   
\[F_t(b\cdot r_\nu:\nu\in\pos(t))=b\cdot F_t(r_\nu:\nu\in\pos(t)),\]
\item[$(\delta)$] if $\langle r_\nu:\nu\in \pos(t)\rangle\subseteq [0,1]$,
$\varepsilon>0$, then there are $r'_\nu>r_\nu$ (for $\nu\in\pos(t)$) 
such that for each $\langle r_\nu'':\nu\in\pos(t)\rangle\subseteq [0,1]$
satisfying $r_\nu\leq r''_\nu< r_\nu'$ (for $\nu\in\pos(t)$) we have
\[F_t(r_\nu'':\nu\in\pos(t))<F_t(r_\nu:\nu\in\pos(t))+\varepsilon.\] 
\end{enumerate}
\end{definition}

From now on (till the end of this section), let $(K,\Sigma,\bF)$ be a fixed
strongly finitary and nice measured tree creating triple for $\bH$. Note
that then the condition (c)$_4$ of Definition \ref{treeforc} is equivalent
to
\begin{enumerate}
\item[(c)$_5$] \quad $(\forall k\in\omega)(\exists n)(\forall\eta\in
T^p)(\lh(\eta)\geq n\ \Rightarrow\ \nor[t_\eta]\geq k)$.
\end{enumerate}

\begin{proposition}
\label{easyprop}
Let $t\in K$. Then: 
\begin{enumerate}
\item[$(\varepsilon)$] If $r_\nu=0$ for $\nu\in\pos(t)$, then $F_t(r_\nu:\nu\in
\pos(t))=0$.
\item[$(\zeta)$] If $\langle r_\nu:\nu\in \pos(t)\rangle\subseteq
[0,1]$, $\varepsilon>0$, then there are $r'_\nu<r_\nu$ (for $\nu\in\pos(t)$) 
such that for each $\langle r_\nu'':\nu\in\pos(t)\rangle\subseteq [0,1]$
satisfying $r_\nu'<r''_\nu\leq r_\nu$ (for $\nu\in\pos(t)$) we have
\[F_t(r_\nu:\nu\in\pos(t))-\varepsilon<F_t(r_\nu'':\nu\in\pos(t)).\]
\end{enumerate}
\end{proposition}

\begin{proof}
$(\varepsilon)$\quad Follows from \ref{nicemix}$(\gamma)$ (take $b=0$).

\noindent $(\zeta)$\quad If $F_t(r_\nu:\nu\in\pos(t))\leq\varepsilon$, then 
any $r_\nu'<r_\nu$ (for $\nu\in\pos(t)$) work. So assume $F_t(r_\nu:\nu\in
\pos(t))>\varepsilon$ and let $b=\frac{F_t(r_\nu:\nu\in\pos(t))-
\varepsilon/2}{F_t(r_\nu:\nu\in\pos(t))}$. Then $0<b<1$. For $\nu\in\pos(t)$
put
\[r_\nu'=\left\{\begin{array}{ll}
-1          &\mbox{ if }r_\nu=0,\\
b\cdot r_\nu&\mbox{ otherwise.}
		\end{array}\right.\]
We are going to show that these $r_\nu'$'s are as required. To this end
suppose that $\langle r''_\nu:\nu\in\pos(t)\rangle\subseteq [0,1]$ is such
that $r'_\nu<r''_\nu\leq r_\nu$ (for all $\nu\in\pos(t)$). Then also $b\cdot
r_\nu\leq r''_\nu$ (for $\nu\in\pos(t)$) and by \ref{nicemix}$(\alpha,
\gamma)$ we get
\[\begin{array}{l}
F_t(r''_\nu:\nu\in\pos(t))\geq F_t(b\cdot r_\nu:\nu\in\pos(t))= b\cdot
F_t(r_\nu:\nu\in\pos(t))=\\
F_t(r_\nu:\nu\in\pos(t))-\varepsilon/2> F_t(r_\nu:\nu\in\pos(t))-
\varepsilon.
  \end{array}\]
\end{proof}

\begin{definition}
\label{mixfor}
Let $p=\langle t^p_\eta:\eta\in T^p\rangle\in\bQ^\tree_\emptyset(K,\Sigma)$.
\begin{enumerate}
\item For a front $A\subseteq T^p$ of $T^p$, we let $T[p,A]=\{\eta\in
T^p:(\exists\rho\in A)(\eta\trianglelefteq\rho)\}$.
\item Let $A$ be a front of $T^p$ and let $f:A\longrightarrow [0,1]$. 
By downward induction on $\eta\in T[p,A]$ we define a mapping 
$\mu^f_{p,A}:T[p,A]\longrightarrow [0,1]$ as follows: 
\begin{itemize}
\item if $\eta\in A$ then $\mu^f_{p,A}(\eta)=f(\eta)$,
\item if $\mu^f_{p,A}(\nu)$ has been defined for all $\nu\in\pos(t^p_\eta)$, 
$\eta\in T[p,A]\setminus A$, then we put $\mu_{p,A}^f(\eta)=F_{t^p_\eta}(
\mu_{p,A}^f(\nu):\nu\in\pos(t^p_\eta))$.
\end{itemize}
\item For $\eta\in T^p$ we define
\[\mu^\bF_p(\eta)=\inf\{\mu^f_{p^{[\eta]},A}(\eta):A\mbox{ is a front of }
(T^p)^{[\eta]}\mbox{ and }f=1_A\},\]
and we let $\mu^\bF(p)=\mu^\bF_p(\mrot(p))$.
\item For $e\in\{\emptyset,4\}$ we let 
\[\bQ^\mtree_e(K,\Sigma,\bF)=\{p\in\bQ^\tree_e(K,\Sigma):\mu^\bF(p)>0\}.\]
It is equipped with the partial order inherited from
$\bQ^\tree_e(K,\Sigma)$.  
\end{enumerate}
\end{definition}

\begin{proposition}
\label{observ}
Assume $p\in\bQ^\tree_\emptyset(K,\Sigma)$ and $A$ is a front of $T^p$.
\begin{enumerate}
\item If $f_0,f_1:A\longrightarrow [0,1]$ are such that $f_0(\nu)\leq
f_1(\nu)$ for all $\nu\in A$, then 
\[(\forall\eta\in T[p,A])(\mu^{f_0}_{p,A}(\eta)\leq\mu^{f_1}_{p,A}
(\eta)).\]
\item If $f_0:A\longrightarrow [0,1]$, $b\in [0,1]$, and $f_1(\nu)=b\cdot
f_0(\nu)$ (for $\nu\in A$), then 
\[(\forall\eta\in T[p,A])(\mu^{f_1}_{p,A}(\eta)=b\cdot\mu^{f_0}_{p,A}
(\eta)).\]
\item If $A'$ is a front of $T^p$ above $A$ (that is, $(\forall\nu'\in A')
(\exists\nu\in A)(\nu\vartriangleleft\nu')$) and $\eta\in T[p,A]$, then
$\mu^{1_{A'}}_{p,A'}(\eta)\leq \mu^{1_A}_{p,A}(\eta)$.
\end{enumerate}
\end{proposition}

\begin{definition}
\label{measure}
Let $p\in\bQ^\tree_\emptyset(K,\Sigma,\bF)$.
\begin{enumerate}
\item A function $\mu:T^p\longrightarrow [0,1]$ is {\em a
semi--$\bF$--measure on $p$} if    
\[(\forall\eta\in T^p)\big(\mu(\eta)\leq F_{t^p_\eta}(\mu(\nu):\nu\in\pos( 
t^p_\eta))\big).\]
\item If above the equality holds for each $\eta\in T^p$, then $\mu$ is
called {\em an $\bF$--measure}. 
\end{enumerate}
\end{definition}

\begin{proposition}
\label{usesemi}
Let $p\in\bQ^\tree_\emptyset(K,\Sigma)$.
\begin{enumerate}
\item If $\mu:T^p\longrightarrow [0,1]$ is a semi--$\bF$--measure on $p$,
then for each $\eta\in T^p$ we have $\mu(\eta)\leq\mu^\bF_p(\eta)$. 
\item If there is a semi--$\bF$--measure $\mu$ on $p$ such that $\mu(
\mrot(p))>0$, then $p\in\bQ^\mtree_\emptyset(K,\Sigma,\bF)$.
\item If $p\in\bQ^\mtree_\emptyset(K,\Sigma,\bF)$, then the mapping
$\eta\mapsto\mu^\bF_p(\eta):T^p\longrightarrow [0,1]$ is an $\bF$--measure
on $p$.  
\end{enumerate}
\end{proposition}

\begin{lemma}
\label{getstr}
Assume $p\in\bQ^\mtree_\emptyset(K,\Sigma,\bF)$ and $0<\varepsilon<1$. Then
there is $\eta\in T^p$ such that $\mu^\bF_p(\eta)\geq 1-\varepsilon$.  
\end{lemma}

\begin{proof}
Assume toward contradiction that $\mu^\bF_p(\eta)<1-\varepsilon$ for all
$\eta\in T^p$. Choose inductively fronts $A_k$ of $T^p$ such that
\begin{itemize}
\item $A_0=\{\mrot(p)\}$,
\item $(\forall\eta\in A_{k+1})(\exists\nu\in A_k)(\nu\vartriangleleft
\eta)$,
\item $\mu^{1_{A_{k+1}}}_{p,A_{k+1}}(\nu)<1-\varepsilon$ for all $\nu\in
A_k$.
\end{itemize}
Note that then (by \ref{observ}(1,2)) for each $k<\omega$ we have
\[\mu^{1_{A_{k+1}}}_{p,A_{k+1}}(\mrot(p))\leq (1-\varepsilon)^{k+1}.\]
Since the right hand side of the inequality above approaches $0$ (as $k\to
\infty$), we get an immediate contradiction with the demand $\mu^\bF(p)>0$. 
\end{proof}

\begin{definition}
\label{normal}
A condition $p\in\bQ^\mtree_\emptyset(K,\Sigma,\bF)$ is called {\em normal}
if for every $\eta\in T^p$ we have $\mu^\bF_p(\eta)>0$. We say that $p$ is
{\em special\/} if for every $\eta\in T^p$ we have $\mu^\bF_p(\eta)\geq
2^{-2^{\lh(\eta)+1}}$. 
\end{definition}

\begin{proposition}
\label{nordense}
\begin{enumerate}
\item Special conditions are dense in $\bQ^\mtree_4(K,\Sigma,\bF)$. (So also
normal conditions are dense.)
\item If $p$ is normal, and $A$ is a front of $T^p$, then $\mu^\bF(p)=
\mu^f_{p,A}(\mrot(p))$, where $f(\nu)=\mu^\bF_p(\nu)$ (for $\nu\in A$).
\end{enumerate}
\end{proposition}

\begin{proof}
1)\quad Let $p\in\bQ^\mtree_4(K,\Sigma,\bF)$; clearly we may assume that
$\nor[t^p_\eta]>1$ for all $\eta\in T^p$. Also we may assume that $\mu^\bF(
p)>3/4$ (remember \ref{getstr}) and $\lh(\mrot(p))>4$.   

Fix  $\eta\in T^p$ such that $\mu^\bF_p(\eta)\geq 2^{-2^{\lh(\eta)}}$ for a
moment. Let $1<a<2$. For each $\nu\in\pos(t^p_\eta)$ pick a front $A_\nu$ of
$(T^p)^{[\nu]}$ such that
\begin{itemize}
\item if $\mu^\bF_p(\nu)<2^{-2^{\lh(\eta)+1}}$, then $\mu^{1_{A_\nu}}_{p^{
[\nu]},A_\nu}(\nu)<2^{-2^{\lh(\eta)+1}}$,
\item if $\mu^\bF_p(\nu)\geq 2^{-2^{\lh(\eta)+1}}$, then $\mu^{1_{A_\nu}}_{
p^{[\nu]},A_\nu}(\nu)<a\mu^\bF_p(\nu)$.
\end{itemize}
Let $X_0=\{\nu\in\pos(t^p_\eta):\mu^\bF_p(\nu)<2^{-2^{\lh(\eta)+1}}\}$,
$X_1=\pos(t^p_\eta)\setminus X_0$, $r_\nu=\mu^{1_{A_\nu}}_{p^{[\nu]},A_\nu}(
\nu)$, and
\[r^\ell_\nu=\left\{\begin{array}{ll}
r_\nu&\mbox{if }\nu\in X_\ell,\\
0    &\mbox{if }\nu\in X_{1-\ell}.\end{array}\right.\]
Apply \ref{nicemix}$(\beta)$ for $t^p_\eta, r^0_\nu,r^1_\nu,r_\nu$ (note
that $F_{t^p_\eta}(r_\nu:r_\nu\in\pos(t^p_\eta))\geq \mu^\bF_p(\eta)\geq
2^{-2^{\lh(\eta)}}$) to pick $s_0^a,s_1^a\in\Sigma(t^p_\eta)$ and $c^a_0,
c^a_1$ such that 
\[c^a_0+c^a_1=(1-2^{-2^{\lh(\eta)}}) F_{t^p_\eta}(r_\nu:\nu\in\pos(
t^p_\eta)),\]
and
\begin{enumerate}
\item[$(\otimes)^a$] if $\ell<2$, $c^a_\ell>0$, then $\nor[s_\ell^a]\geq
\nor[t^p_\eta]-1$, $\pos(s_\ell^a)\subseteq X_\ell$, and 
\[F_{s_\ell^a}(r_\nu:\nu\in\pos(s_\ell^a))\geq c^a_\ell.\]
\end{enumerate}
Note that, if $c^a_0>0$, then $c^a_0\leq F_{s_0^a}(r_\nu:\nu\in\pos(
s_\ell^a))\leq 2^{-2^{\lh(\eta)+1}}$, and thus
\[c^a_1\geq (1-2^{-2^{\lh(\eta)}}) F_{t^p_\eta}(r_\nu:\nu\in
\pos(t^p_\eta))-2^{-2^{\lh(\eta)+1}}>0.\]
Also, letting $r^*_\nu=\min\{a\mu^\bF_p(\nu),1\}$, 
\[F_{s_1^a}(r_\nu:\nu\in\pos(s^a_1))\leq F_{s_1^a}(r^*_\nu:\nu\in
\pos(s^a_1))\leq a\cdot F_{s_1^a}(\mu^\bF_p(\nu):\nu\in\pos(s^a_1)).\]
Together
\begin{enumerate}
\item[$(*)_a$] $(1-2^{-2^{\lh(\eta)}}) F_{t^p_\eta}(r_\nu:\nu\in
\pos(t^p_\eta))-2^{-2^{\lh(\eta)+1}}\leq a\cdot F_{s_1^a}(\mu^\bF_p(\nu):\nu\in\pos(
s^a_1))$.
\end{enumerate}
Since $(K,\Sigma)$ is strongly finitary, considering $a\to 1$ (and using
\ref{nicemix}$(\delta)$), we find $s_\eta\in\Sigma(t_\eta^p)$ such that
$\nor[s_\eta]\geq\nor[t^p_\eta]-1$ and $\mu^\bF_p(\nu)\geq 2^{-2^{\lh(\eta)+
1}}$ for all $\nu\in\pos(s_\eta)$, and
\[\mu^\bF_p(\eta)=F_{t^p_\eta}(\mu^\bF_p(\nu):\nu\in\pos(t^p_\eta))\leq
\frac{F_{s_\eta}(\mu^\bF_p(\nu):\nu\in\pos(s_\eta))+2^{-2^{\lh(\eta)+1}}}{1-
2^{-2^{\lh(\eta)}}}.\]
Note that also, as $2^{-2^{\lh(\eta)}}\leq\mu^\bF_p(\eta)$, 
\[\mu^\bF_p(\eta)(1-2^{-2^{\lh(\eta)}})-2^{-2^{\lh(\eta)+1}}\geq 
\mu^\bF_p(\eta)(1-2^{1-2^{\lh(\eta)}}),\]
so
\begin{enumerate}
\item[$(**)$] \quad $\mu^\bF_p(\eta)\cdot (1-2^{1-2^{\lh(\eta)}})\leq
F_{s_\eta}(\mu^\bF_p(\nu):\nu\in\pos(s_\eta))$. 
\end{enumerate}
\medskip

Now, starting with $\mrot(p)$, build a tree $S$ and a system $q=\langle
s_\eta:\eta\in S\rangle$ such that $\suc_S(\eta)=\pos(s_\eta)$. It should be
clear that in this way we will get a condition in $\bQ^\tree_4(K,\Sigma)$
(stronger than $p$).  Why is $q$ in $\bQ^\mtree_4(K,\Sigma,\bF)$? Let
$k^*>\lh(\mrot(q))$, $A=\{\nu\in S:\lh(\nu)=k^*\}$ and $f=1_A$. Using
$(**)$, we may show by downward induction that for every $\eta\in T[q,A]$ we
have 
\[\begin{array}{r}
\mu^f_{q,A}(\eta)\geq\mu^\bF_p(\eta)\cdot\prod\limits_{k=\lh(\eta)}^{k^*
-1} (1-2^{1-2^k})\geq\mu^\bF_p(\eta)\cdot (1-2^{2-2^{\lh(\eta)}})\geq {\ }\\
 2^{-2^{\lh(\eta)}}(1-2^{2-2^{\lh(\eta)}})\geq 2^{-2^{\lh(\eta)+1}}.
  \end{array}\]
Now we may easily conclude that $q\in\bQ^\mtree_4(K,\Sigma,\bF)$ is
special. 
\medskip

\noindent 2)\quad Let $A$ be a front of $T^p$, $p$ normal (so, in
particular, $\mu^\bF_p(\nu)>0$ for $\nu\in A$). Fix $a>1$ for a moment.

For each $\nu\in A$ pick a front $A_\nu$ of $(T^p)^{[\nu]}$ such that
$\mu^{1_{A_\nu}}_{p,A_\nu}(\nu)<a\cdot\mu^\bF_p(\nu)$. Let
$B=\bigcup\limits_{\nu\in A} A_\nu$ and $f(\nu)=\mu^\bF_p(\nu)$ for $\nu\in
A$. By downward induction one can show that for all $\rho\in T[p,A]$ we have
$\mu^{1_B}_{p,B}(\rho)\leq a\cdot \mu^f_{p,A}(\rho)$. Then, in particular,
we have 
\[\mu^\bF(p)\leq\mu^{1_B}_{p,B}(\mrot(p))\leq a\cdot\mu^f_{p,A}(\mrot(p)),\]
and hence (letting $a\to 1$) $\mu^\bF(p)\leq\mu^f_{p,A}(\mrot(p))$. The
reverse inequality is even easier (remember \ref{observ}(1)).
\end{proof}

\begin{lemma}
\label{useful}
Let $p\in\bQ^\mtree_4(K,\Sigma)$ be a normal condition such that $\mu^\bF(p)
>\frac{1}{2}$, $\nor[t^p_\eta]>2$ for all $\eta\in T^p$, and let $k_0=\lh(
\mrot(p))>4$, $0<\varepsilon\leq 2^{-(1+k_0)}$. Suppose that $B$ is an
antichain of $T^p$, and that for each $\nu\in B$ we are given a normal
condition $q_\nu\geq p^{[\nu]}$ such that 
\[\mrot(q_\nu)=\nu\quad\mbox{ and }\quad \mu^\bF(q_\nu)\geq 1-\varepsilon.\] 
Then at least one of the following conditions holds. 
\begin{enumerate}
\item[(i)] There is a normal condition $q\in\bQ^\mtree_4(K,\Sigma,\bF)$ such 
that 
\[q\geq p,\quad \mrot(q)=\mrot(p),\quad\mbox{ and  }\quad T^q\cap B=
\emptyset.\]
\item[(ii)]  There is a normal condition $q\in\bQ^\mtree_4(K,\Sigma,\bF)$
such that 
\begin{itemize}
\item $q\geq p$, $\mrot(q)=\mrot(p)$, $\mu^\bF(q)\geq (1-2^{-k_0})
\mu^\bF(p)$, and 
\item $T^q\cap B$ is a front of $T^q$, and $q^{[\nu]}=q_\nu$ for $\nu\in
T^q\cap B$, and
\item if $\eta\in T^q$, $\eta\vartriangleleft\nu\in B$, then $\nor[t^q_\eta] 
\geq\nor[t^p_\eta]-2$.  
\end{itemize}
\end{enumerate}
\end{lemma}

\begin{proof}
Let $e_\ell=2^{1-2^\ell}$ (for $\ell<\omega$); note that $(e_\ell)^2=2
e_{\ell+1}$.  

Fix $k>\lh(\mrot(p))$ for a while. Let $A$ be a front of $T^p$ such that
\[\{\nu\in B:\lh(\nu)\leq k\}\subseteq A\quad\mbox{and}\quad (\forall\nu\in 
A)(\nu\notin B\ \Rightarrow\ \lh(\nu)=k).\]
By downward induction, for each $\nu\in T[p,A]$, we define $r_\nu^0,r^1_\nu
\in [0,1]$ and $s^0_\nu,s^1_\nu\in\Sigma(t^p_\nu)$ such that
\begin{enumerate}
\item[$(\alpha)$] If $\nu\in A\cap B$, then $r^0_\nu=0$, $r_\nu^1=\mu^\bF(
q_\nu)$. 
\item[$(\beta)$] If $\nu\in A\setminus B$, then $r_\nu^0=\mu^\bF_p(\nu)$,
$r^1_\nu=0$. 
\item[$(\gamma)$] If $\nu\in T[p,A]\setminus A$, $\lh(\nu)=m$, then:

if $\mu^\bF_p(\nu)\cdot (1-\varepsilon)\cdot\prod\limits_{\ell=m}^{k-1}(1-3 
e_\ell)<e_m$, then $r^0_\nu=r^1_\nu=0$,

else $r^0_\nu+r^1_\nu\geq\mu^\bF_p(\nu)\cdot (1-\varepsilon)\cdot
\prod\limits_{\ell=m}^{k-1}(1-3e_\ell)$.
\end{enumerate}
Clauses $(\alpha),(\beta)$ define $r^0_\nu,r^1_\nu$ for $\nu\in A$;
$s^0_\nu,s^1_\nu$ are not defined then (or are arbitrary). 

Suppose $\eta\in T[p,A]\setminus A$, $\lh(\eta)=k-1$. If $\mu^\bF_p(\eta)
\cdot (1-\varepsilon)\cdot (1-3e_{k-1})<e_{k-1}$, then we let $r^0_\eta=
r^1_\eta=0$ (and $s^0_\eta,s^1_\eta$ are not defined). So assume now that  
\[\mu^\bF_p(\eta)\cdot (1-\varepsilon)\cdot (1-3e_{k-1})\geq e_{k-1}.\]
Then also (as $r^0_\nu+r^1_\nu\geq\mu^\bF_p(\nu)\cdot (1-\varepsilon)$)
\[\begin{array}{rl}
F_{t^p_\eta}(r^0_\nu+r^1_\nu:\nu\in\pos(t^p_\eta))\geq\mu^\bF_p(
\eta)\cdot (1-\varepsilon)\geq&\ \\
 \mu^\bF_p(\eta)\cdot (1-\varepsilon)\cdot(1-3e_{k-1})\geq& e_{k-1}>
2^{-2^{k-1}},\end{array}\]
and we may apply \ref{nicemix}$(\beta)$ to pick $r^0_\eta,r^1_\eta$ and
$s^0_\eta,s^1_\eta\in\Sigma(t^p_\eta)$ such that
\begin{enumerate}
\item[(i)] $r^0_\eta+r^1_\eta\geq (1-e_{k-1})\cdot F_{t^p_\eta}(r^0_\nu+
r^1_\nu:\nu\in\pos(t^p_\eta))\geq\mu^\bF_p(\eta)\cdot (1-\varepsilon)\cdot
(1-3e_{k-1})$, 
\item[(ii)] if $r^\ell_\eta>0$, $\ell<2$, then $\nor[s^\ell_\eta]\geq
\nor[t^p_\eta]-1$, $\pos(s^\ell_\eta)\subseteq\{\nu\in\pos(t^p_\eta):
r^\ell_\nu>0\}$, and $F_{s^\ell_\eta}(r^\ell_\nu:\nu\in\pos(s^\ell_\eta))
\geq r^\ell_\eta$.
\end{enumerate}

Suppose now that $\eta\in T[p,A]\setminus A$, $\lh(\eta)=m-1<k-1$, and
$r^0_\nu,r^1_\nu$ have been defined for all $\nu\in\pos(t^p_\eta)$ (and they
satisfy clause $(\gamma)$). If 
\[\mu^\bF_p(\eta)\cdot (1-\varepsilon)\cdot\prod\limits_{\ell=m-1}^{k-1}(1-
3 e_\ell) <e_{m-1},\]
then we let $r^0_\eta=r^1_\eta=0$ (and $s^0_\eta,s^1_\eta$ are not
defined). So assume
\[\mu^\bF_p(\eta)\cdot (1-\varepsilon)\cdot\prod\limits_{\ell=m-1}^{k-1}(1-3 
e_\ell) \geq e_{m-1}.\] 
Then for $\nu\in\pos(t^p_\eta)$ we let
\[r^*_\nu=\left\{\begin{array}{ll}
r^0_\nu+r^1_\nu&\mbox{ if }r^0_\nu+r^1_\nu>0,\\
e_m&\mbox{ otherwise,}\end{array}\right.\]
and we note that 
\[F_{t^p_\eta}(r^*_\nu:\nu\in\pos(t^p_\eta))\geq\mu^\bF_p(\eta)\cdot
(1-\varepsilon)\cdot \prod\limits_{\ell=m}^{k-1}(1-3e_\ell)\geq e_{m-1}>
2^{-2^{m-1}}.\] 
Applying \ref{nicemix}$(\beta)$ choose $t^0,t^1\in\Sigma(t^p_\eta)$ and
$c_0,c_1$ such that $c_0+c_1\geq (1-e_{m-1}) F_{t^p_\eta}(r^*_\nu:\nu\in
\pos(t^p_\eta))$ and
\begin{itemize}
\item if $c_0>0$, then $\pos(t^0)\subseteq\{\nu\in\pos(t^p_\eta):r^0_\nu+
r^1_\nu=0\}$, $\nor[t^0]\geq \nor[t^p_\eta]-1$ and $F_{t^0}(r^*_\nu:\nu\in
\pos(t^0))\geq c_0$, 
\item  if $c_1>0$, then $\pos(t^1)\subseteq\{\nu\in\pos(t^p_\eta):r^0_\nu+
r^1_\nu>0\}$, $\nor[t^1]\geq \nor[t^p_\eta]-1$ and $F_{t^1}(r^*_\nu:\nu\in
\pos(t^1))\geq c_1$.
\end{itemize}
Now look at the definition of $r^*_\nu$. If $c_0>0$, then $F_{t^0}(r^*_\nu:
\nu\in\pos(t^0))\leq e_m$, so $c_0\leq e_m\leq (e_{m-1})^2$. Therefore 
\[\begin{array}{l}
F_{t^1}(r^*_\nu:\nu\in\pos(t^1))\geq c_1\geq (1-e_{m-1})\cdot\mu^\bF_p(\eta)
\cdot (1-\varepsilon)\cdot\prod\limits_{\ell=m}^{k-1}(1-3e_\ell) -e_m\geq\\
(1-e_{m-1})\mu^\bF_p(\eta)(1-\varepsilon)\cdot\prod\limits_{\ell=m}^{k-1}
(1-3e_\ell)-e_{m-1}\mu^\bF_p(\eta)(1-\varepsilon)\cdot \prod\limits_{\ell=
m}^{k-1}(1- 3e_\ell)=\\
\mu^\bF_p(\eta)(1-\varepsilon)\cdot\prod\limits_{\ell=m}^{k-1}(1-3e_\ell)
\cdot (1-2e_{m-1})\geq e_{m-1}\cdot (1-2e_{m-1})>2^{-2^{m-1}}.
\end{array}\]
Hence we may apply \ref{nicemix}$(\beta)$ again and get $r^0_\eta,r^1_\eta$
and $s^0_\eta,s^1_\eta\in\Sigma(t^1)\subseteq\Sigma(t^p_\eta)$ such that 
\[\begin{array}{l}
r^0_\eta+r^1_\eta\geq (1-e_{m-1})\cdot F_{t^1}(r^*_\nu:\nu\in\pos(t^1))\geq\\
\mu^\bF_p(\eta)(1-\varepsilon)\cdot\prod\limits_{\ell=m}^{k-1}(1-3e_\ell)\cdot
(1-2e_{m-1})\cdot (1-e_{m-1})\geq\\
\mu^\bF_p(\eta)\cdot (1-\varepsilon)\cdot\prod\limits_{\ell=m-1}^{k-1}(1
-3e_\ell),\end{array}\]
and if $r^\ell_\eta>0$, $\ell<2$, then $\pos(s^\ell_\eta)\subseteq\{\nu\in
\pos(t^p_\eta):r^\ell_\nu>0\}$, $\nor[s^\ell_\eta]\geq \nor[t^p_\eta]-2$ and
$F_{s^\ell_\eta}(r^\ell_\nu:\nu\in\pos(s^\ell_\eta))\geq r^\ell_\eta$. 

Note that (as $k_0>4$) 
\[\varepsilon+\sum\limits_{\ell=k_0}^{k-1}3e_\ell\leq\frac{1}{2^{k_0+1}}+6
\cdot\sum\limits_{\ell=k_0}^\infty\frac{1}{2^{2^\ell}}\leq
 \frac{3}{2^{k_0+2}}.\] 
 Therefore,
\[\begin{array}{r}
\mu^\bF_p(\mrot(p))\cdot\prod\limits_{\ell=k_0}^{k-1}(1-3e_\ell)
\cdot (1-\varepsilon)\geq \mu^\bF_p(\mrot(p))\cdot (1-(\varepsilon+
\sum\limits_{\ell=k_0}^{k-1}3e_\ell))\geq\\
\mu^\bF_p(\mrot(p))\cdot(1-\frac{3}{2^{k_0+2}})>\frac{1}{2}\cdot
\frac{29}{32}>e_{k_0}.
  \end{array}\] 
Hence also (by $(\gamma)$)
\[\mu^\bF_p(\mrot(p))(1-\frac{3}{2^{k_0+2}})\leq \mu^\bF_p(\mrot(p))
\cdot\prod\limits_{\ell=k_0}^{k-1}(1-3e_\ell)\cdot (1-\varepsilon)\leq 
r^0_{\mrot(p)}+r^1_{\mrot(p)}.\]
Now, if $r^\ell_{\mrot(p)}>0$, $\ell<2$, then we build inductively a finite
tree $S^k_\ell\subseteq T[p,A]$ as follows. We declare that $\mrot(
S^k_\ell)=\mrot(p)$, $s^{\ell,k}_{\mrot(p)}=s^\ell_{\mrot(p)}$, and
$\suc_{S^k_\ell}(\mrot(p))=\pos(s^{\ell,k}_{\mrot(p)})$. If we have decided
that $\eta\in S^k_\ell$, $\eta\notin A$ (and $r^\ell_\eta>0$), then we also
declare $s^{\ell,k}_\eta=s^\ell_\eta$, $\suc_{S^k_\ell}(\eta)=\pos(s^{\ell,
k}_\eta)$ (note $r^\ell_\nu>0$ for $\nu\in\pos(s^{\ell,k}_\eta)$).

Then, if $S^k_0$ is defined, $S^k_0\cap B=\emptyset$, and, if $S^k_1$ is
defined, $S^k_1\cap A\subseteq B$. Also, if we ``extend'' $S^k_0$ using
$p^{[\nu]}$ (for $\nu\in S^k_0\cap A$), then we get a condition $q^k_0\geq
p$ such that $\mu^\bF(q^k_0)\geq r^0_{\mrot(p)}\stackrel{\rm def}{=}
r^{0,k}$. Likewise, if we ``extend'' $S^k_1$ using $q_\nu$ (for $\nu\in
S^k_1\cap A$), then we get a condition $q^k_1\geq p$ such that $\mu^\bF(
q^k_1)\geq r^1_{\mrot(p)}\stackrel{\rm def}{=}r^{1,k}$. 
\medskip

If for some $k>\lh(\mrot(p))$ we have $r^{1,k}\geq(1-2^{-k_0})\mu^\bF(p)$,
then we use the respective condition $q^k_1$ to witness the demand (ii) of
the lemma. So assume that for each $k>\lh(\mrot(p))$ we have $r^{1,k}<(1- 
2^{-k_0})\mu^\bF(p)$, and thus
\[r^{0,k}>(1-\frac{3}{2^{k_0+2}})\mu^\bF(p)-(1-\frac{1}{2^{k_0}})\mu^\bF(p)
=\frac{1}{2^{k_0}}\mu^\bF(p)>0.\] 
Apply the K\"onig Lemma to find an infinite set $I\subseteq\omega\setminus 
(k_0+1)$ such that for all $k,k',k''\in I$, $k<k'<k''$, we have
\[(\forall\eta\in S^{k'}_0)(\lh(\eta)\leq k\ \Rightarrow\ \eta\in S^{k''}_0\
\&\ s^{0,k'}_\eta=s^{0,k''}_\eta).\]
Then $S^q=\{\eta:(\forall^\infty k\in I)(\eta\in S^k_0)\}$, $s^q_\eta=
s^{0,k}_\eta$ (for sufficiently large $k\in I$) determine a condition $q$
witnessing the first assertion of the lemma.  
\end{proof}

\begin{lemma}
\label{preconrea}
Assume that $\dot{\tau}$ is a $\bQ^\mtree_4(K,\Sigma,\bF)$--name for an
ordinal, $n\leq m<\omega$ and $p\in\bQ^\mtree_4(K,\Sigma,\bF)$ is a normal 
condition such that $\mu^\bF(p)>\frac{1}{2}$, and $\nor[t^p_\eta]>n+2$ for 
$\eta\in T^p$. Let $k_0=\lh(\mrot(p))>2$. Then there is a normal condition
$q\in\bQ^\mtree_4(K,\Sigma,\bF)$ such that  
\begin{enumerate}
\item[(a)] $q\geq p$, $\mrot(q)=\mrot(p)$, $\mu^\bF(q)\geq (1-2^{-k_0})
\mu^\bF(p)$, and 
\item[(b)] $(\forall\eta\in T^q)(\nor[t^q_\eta]\geq n)$, and 
\item[(c)] there is a front $A$ of $T^q$ such that for every $\nu\in A$:
\begin{itemize}
\item the condition $q^{[\nu]}$ forces a value to $\dot{\tau}$,
\item $\mu^\bF_q(\nu)>\frac{1}{2}$, $\lh(\nu)>k_0$,
\item if $\nu\trianglelefteq\eta\in T^q$, then $\nor[t^q_\eta]\geq m$.
\end{itemize}
\end{enumerate}
\end{lemma}

\begin{proof}
Let $B$ consist of all $\nu\in T^p$ such that 
\begin{enumerate}
\item[$(\alpha)$] $\lh(\nu)>k_0$ and there is a normal condition
$q\in\bQ^\mtree_4(K,\Sigma,\bF)$ stronger than $p^{[\nu]}$ and such that
$\mrot(q)=\nu$, $\mu^\bF(q)\geq (1-2^{-(2+k_0)})$, $(\forall\eta\in
T^q)(\nor[t^q_\eta]\geq m)$, and for some front $A$ of $T^q$, for every
$\eta\in A$: 
\begin{enumerate}
\item[$(\otimes)$] \qquad $\mu^\bF_q(\eta)>1/2$ and the condition
$q^{[\eta]}$ decides the value of $\dot{\tau}$,
\end{enumerate}
and 
\item[$(\beta)$] no initial segment of $\eta$ has the property stated in
$(\alpha)$ above. 
\end{enumerate}
Note that $B$ is an antichain of $T^p$, and $B\cap T^{p'}\neq\emptyset$ for
every condition $p'\geq p$ such that $\mrot(p')=\mrot(p)$ (by
\ref{getstr}). For each $\nu\in B$ fix a condition $q_\nu$ witnessing clause
$(\alpha)$ (for $\nu$). Now apply \ref{useful}: case (i) there is not
possible by what we stated above, so we get a condition $q$ as described in
\ref{useful}(ii). It should be clear that it is as required here. 
\end{proof}

\begin{theorem}
\label{contreading}
Suppose that $p\in\bQ^\mtree_4(K,\Sigma,\bF)$, and $\dot{\tau}_n$ are
$\bQ^\mtree_4(K,\Sigma,\bF)$--names for ordinals ($n<\omega$). Then there
are a condition $q\geq p$ and fronts $A_n$ of $T^q$ (for $n<\omega$) such
that for each $n<\omega$ and $\nu\in A_n$, the condition $q^{[\nu]}$ decides
the value of $\dot{\tau}_n$. 
\end{theorem}

\begin{proof}
We may assume that $p$ is normal, $k_0=\lh(\mrot(p))>2$, $\mu^\bF(p)>
\frac{1}{2}$, and $\nor[t^p_\eta]>3$ for $\eta\in T^p$. We build inductively
a sequence $\langle q_n,A_n:n<\omega\rangle$ such that   
\begin{enumerate}
\item $q_n\in\bQ^\mtree_4(K,\Sigma,\bF)$ is a normal condition, $\mrot(q_n)
=\mrot(p)$, $q_n\leq q_{n+1}$, $q_0=p$,
\item $A_n\subseteq T^{q_{n+1}}$ is a front of $T^{q_{n+1}}$, $(\forall\nu
\in A_n)(\exists \eta\in A_{n+1})(\nu\vartriangleleft\eta)$, 
\item if $\nu\in A_n$, then $\mu^\bF_{q_{n+1}}(\nu)>\frac{1}{2}$, and for
each $\eta\in T^{q_{n+1}}$ such that $\nu\trianglelefteq\eta$ we have
$\nor[t^{q_{n+1}}_\eta]\geq n+4$, 
\item if $\mrot(p)\trianglelefteq\eta\trianglelefteq\nu\in A_n$, then
$t^{q_{n+1}}_\eta=t^{q_{n+2}}_\eta$,
\item for each $\nu\in A_n$, the condition $(q_{n+1})^{[\nu]}$ decides the
value of $\dot{\tau}_n$,
\item $\mu^\bF(q_{n+1})\geq\prod\limits_{\ell=k_0}^{k_0+n}(1-2^{-\ell})\cdot 
\mu^\bF(p)$. 
\end{enumerate}
The construction can be carried out by \ref{preconrea} ($q_1,A_0$ are
obtained by applying \ref{preconrea} to $p$ and $\dot{\tau}_0$; if
$q_{n+1},A_n$ have been defined, then we apply \ref{preconrea} to
$\dot{\tau}_{n+1}$ and $(q_{n+1})^{[\nu]}$ for $\nu\in A_n$; remember
\ref{observ}). Next define $q=\langle t^q_\eta:\eta\in T^q\rangle$ so that
$\mrot(q)=\mrot(p)$, each $A_n$ is a front of $T^q$, and if $\mrot(p)
\trianglelefteq\eta\trianglelefteq\nu\in A_n$ then $t^q_\eta=t^{q_{n+
1}}_\eta$. It is straightforward to check that $q$ is as required in
\ref{contreading}. 
\end{proof}

\begin{corollary}
\label{corprobou}
Let $(K,\Sigma,\bF)$ be a strongly finitary nice measured tree creating
triple. Then the forcing notion $\bQ^\mtree_4(K,\Sigma,\bF)$ is proper and
$\baire$--bounding. 
\end{corollary}

\begin{proposition}
\label{snep}
If $(K,\Sigma,\bF)$ is a strongly finitary nice measured tree creating
triple, then the forcing notion $\bQ^\mtree_4(K,\Sigma,\bF)$ is explicitly
$\aleph_0$--snep (see \cite[Definitions 1.5, 1.8]{Sh:630}.
\end{proposition}

\begin{proof}
Should be clear at the moment; compare the proof of 
\cite[Lemma 3.1]{Sh:630}. 
\end{proof}

\section{The Forcing}
In this section we define a nice, strongly finitary measured tree creating
triple $(K^*,\Sigma^*,\bF^*)$, and we show several technical properties of
it and of the forcing notion $\starQ$. This forcing will be used in the next
two sections to show our main results \ref{main} and \ref{twomain}.
\medskip

Let $k<\omega$. Fix a function $\varphi_k:\omega\longrightarrow\omega$ such
that 
\[\varphi_k(0)=2^{k+4}\quad\mbox{and}\quad \varphi_k(i+1)>2^{2^{k+3}}+
\varphi_k(i)+\frac{2^{2k+7}}{\log_2(1+2^{-2^{2k+7}})}.\] 
Let $N_k=2^{1+\lfloor\log_2(\varphi_k(k+1))\rfloor}$ (where $\lfloor
r\rfloor$ is the integer part of the real number $r$), and let $\bH^*(k)=
2^{\textstyle N_k}$. 
\medskip 

Let $K^*$ consist of tree creatures $t\in \TCR[\bH^*]$ such that
\begin{itemize}
\item $\dis[t]=(k_t,\eta_t,n_t,g_t,P_t)$, where $n_t\leq k_t<\omega$,
$\eta_t\in\prod\limits_{i<k_t}\bH^*(i)$, $g_t$ is a partial function from
$N_{k_t}$ to $2$ such that $|g_t|\leq\varphi_{k_t}(k_t-n_t)$, and 
\[\emptyset\neq P_t\subseteq\{f\in\bH^*(k_t):g_t\subseteq f\},\]
\item $\nor[t]=n_t$,
\item $\val[t]=\{\langle\eta_t,\nu\rangle:\eta_t\vartriangleleft\nu\in
\prod\limits_{i\leq k_t}\bH^*(i)\ \&\ \nu(k_t)\in P_t\}$.
\end{itemize}
\medskip

The operation $\Sigma^*$ is trivial, and for $t\in K^*$:
\[\Sigma^*(t)=\{s\in K^*:\eta_s=\eta_t\ \&\ n_s\leq n_t\ \&\ g_t\subseteq
g_s\ \&\ P_s\subseteq P_t\}.\]

Finally, for $t\in K^*$ and $\langle r_\nu:\nu\in\pos(t)\rangle\subseteq
[0,1]$ we let 
\[\begin{array}{lr}
F^*_t(r_\nu:\nu\in\pos(t))=&\\
\min\{2^{|h|-N_k}\cdot\sum\{r_\nu: h\subseteq\nu(k_t)\in P_t\}:&
h\mbox{ is a partial function from $N_{k_t}$ to }2,\ \ \\
& g_t\subseteq h\ \mbox{ and }\ |h\setminus g_t|\leq 2^{k_t+3}\}.
  \end{array}\]
(So this defines $\bF^*=\langle F^*_t:t\in K^*\rangle$.)
\medskip

It should be clear that $(K^*,\Sigma^*,\bF^*)$ is a strongly finitary
measured tree creating triple. (And now we are aiming at showing that it is
nice.) 

\begin{lemma}
\label{Sah7.2}
Assume that $t\in K^*$, $\nor[t]>1$, and $g'$ is a partial function from
$N_{k_t}$ to $2$ such that $g'\supseteq g_t$ and $|g'\setminus g_t|\leq
2^{k_t+3}$. Furthermore, suppose that $r_\nu\in[0,1]$ (for $\nu\in\pos(t)$)
are such that
\[2^{-2^{k_t+3}}\leq 2^{|g'|-N_{k_t}}\cdot\sum\{r_\nu:\nu\in\pos(t)\ \&\ g'
\subseteq\nu(k_t)\}\stackrel{\rm def}{=}a.\]
Then there is $s\in\Sigma^*(t)$ such that
\begin{enumerate}
\item[$(\alpha)$] $\nor[s]=\nor[t]-1$, $g'\subseteq g_s$,
\item[$(\beta)$]  $F^*_s(r_\nu:\nu\in\pos(s))\geq a\cdot
(1-2^{-2^{k_t+3}})$, 
\item[$(\gamma)$] if $h$ is a partial function from $N_{k_s}$ to $2$ such
that $g_s\subseteq h$ and $|h\setminus g_s|\leq 2^{k_s+3}$, then 
\[\frac{\sum\{r_\nu:\nu\in\pos(s)\ \&\ h\subseteq\nu(k_s)\}}{2^{N_{k_s}-
|h|}}\]
is in the interval $[F^*_s(r_\nu:\nu\in\pos(s)),F^*_s(r_\nu:\nu\in\pos(s))
\cdot(1+2^{-2^{k+3}})]$.
\end{enumerate}
\end{lemma}

\begin{proof}
Let $k=k_t$, $n=n_t$.

We try to choose inductively partial functions $g_\ell$ from $N_k$ to $2$
such that  
\begin{enumerate}
\item[(a)]  $g'=g_0\subseteq g_1 \subseteq \ldots$, $|g_\ell\setminus g'|
\leq \ell\cdot 2^{k+3}$, 
\item[(b)$_\ell$] $2^{|g_\ell|-N_k}\cdot \sum\{r_\nu:\nu\in\pos(t)\ \&\
g_\ell\subseteq\nu(k)\}\geq a\cdot (1+2^{-2^{2k+7}})^\ell$.
\end{enumerate}
Note that in (b)$_\ell$, the left hand side expression is not more than $1$, 
so if the inequality holds, then (as $a\geq 2^{-2^{k+3}}$)
\begin{enumerate}
\item[$(\oplus)$]\qquad\qquad $\displaystyle \ell\leq \frac{2^{k+3}}{
\log_2(1+2^{-2^{2k+7}})}$.
\end{enumerate}
Consequently, in the procedure described above, we are stuck at some
$\ell_0$ satisfying $(\oplus)$. Let 
\[g_s=g_{\ell_0},\quad n_s=n-1,\quad k_s=k,\quad \eta_s=\eta_t,\quad P_s=\{f
\in P_t:g_{\ell_0}\subseteq f\}.\]
So this defines $s$, but we have to check that $s\in K^*$. For this note
that 
\[|g_s|\leq |g'|+\ell_0\cdot 2^{k+3}\leq \varphi_k(k-n)+2^{k+3}+\frac{2^{2k
+6}}{\log_2(1+2^{-2^{2k+7}})}\leq\varphi_k(k-n_s).\]
(So indeed $s\in K^*$, and plainly $s\in\Sigma^*(t)$.) Also note that 
\[2^{|g_s|-N_k}\cdot \sum\{r_\nu:\nu\in\pos(s)\}\geq a\cdot (1+2^{-2^{2k+7}}
)^{\ell_0}\stackrel{\rm def}{=}a^*\geq a.\]
Now, suppose that $u\subseteq N_k\setminus\dom(g_s)$, $|u|\leq 2^{k+3}$. Let
$h:u\longrightarrow 2$. We cannot use $g_s\conc h$ as $g_{\ell_0+1}$, so the
condition (b)$_{\ell_0+1}$ fails for it. Therefore
\[\begin{array}{r}
b_h\stackrel{\rm def}{=} 2^{|g_s|+|h|-N_k}\cdot\sum\{r_\nu:\nu\in\pos(t)\
\&\ g_s\conc h\subseteq\nu(k)\}<{}\\
\ \\
a\cdot (1+2^{-2^{2k+7}})^{\ell_0+1}=a^*\cdot (1+2^{-2^{2k+7}}).
  \end{array}\]
\begin{claim}
\label{cl1}
For each $h:u\longrightarrow 2$, we have
\[b_h\geq a^*\cdot (1-2^{-2^{k+4}}).\]
\end{claim}

\begin{proof}[Proof of the claim]
Assume that $h_0:u\longrightarrow 2$ is such that $b_{h_0}<a^*\cdot (1-2^{
-2^{k+4}})$. We know that $b_h<a^*\cdot (1+2^{-2^{2k+7}})$ for each $h:u
\longrightarrow 2$, so 
\[\begin{array}{l}
a^*\cdot 2^{N_k-|g_s|}\leq \sum\{r_\nu:\nu\in\pos(s)\}\leq {}\\
a^*\cdot (1-2^{-2^{k+4}})\cdot 2^{N_k-|g_s|-|u|}+ a^*\cdot (1+2^{-2^{2k+7}}) 
\cdot (2^{|u|}-1)\cdot 2^{N_k-|g_s|-|u|}.
  \end{array}\]
Hence
\[\begin{array}{r}
2^{|u|}\leq (1-2^{-2^{k+4}})+(1+2^{-2^{2k+7}})\cdot (2^{|u|}-1)=\quad\\
2^{|u|}(1+2^{-2^{2k+7}})-(2^{-2^{k+4}}+2^{-2^{2k+7}}),
  \end{array}\]
and so $2^{-2^{k+4}}\leq 2^{-2^{k+4}}+2^{-2^{2k+7}}\leq 2^{|u|}\cdot 2^{-
2^{2k+7}}\leq 2^{2^{k+3}-2^{2k+7}}$, a contradiction.
\end{proof}

Consequently, we get that 
\[F^*_s(r_\nu:\nu\in\pos(s))\geq a^*\cdot (1-2^{-2^{k+4}})\geq a\cdot
(1-2^{-2^{k+3}}),\]
so $s$ satisfies the demand $(\beta)$.

But we also know that for each partial function $h$ from $N_k$ to $2$, if
$g_s\subseteq h$ and $|h\setminus g_s|\leq 2^{k+3}$, then 
\[\begin{array}{r}
b_h<a^*\cdot (1+2^{-2^{2k+7}})\leq F^*_s(r_\nu:\nu\in\pos(s))\frac{1+
2^{-2^{2k+7}}}{1-2^{-2^{k+4}}}\leq\quad\\ 
F^*_s(r_\nu:\nu\in\pos(s))\cdot (1+2^{-2^{k+3}}),
\end{array}\]
and thus $s$ satisfies the demand $(\gamma)$ as well.
\end{proof}

\begin{proposition}
\label{nice}
$(K^*,\Sigma^*,\bF^*)$ is a nice (strongly finitary) measured tree creating
triple. 
\end{proposition}

\begin{proof}
Clauses \ref{nicemix}($\alpha,\gamma,\delta$) should be obvious, so let us
check \ref{nicemix}$(\beta)$ only.

Let $t\in K^*$, $k=k_t$, $r^0_\nu,r^1_\nu,r_\nu$ be as in the assumptions of
\ref{nicemix}$(\beta)$. So in particular
\[2^{|g_t|-N_k}\cdot \sum\{r_\nu:\nu\in\pos(t)\}\geq F^*_t(r_\nu:\nu\in
\pos(t))\geq 2^{-2^k}>2^{-2^{k+3}}.\]
For $\ell<2$ let $a_\ell=2^{|g_t|-N_k}\cdot\sum\{r^\ell_\nu:\nu\in\pos(t)
\}$.

First, we consider the case when both $a_0$ and $a_1$ are not smaller than
$2^{-2^{k+3}}$. Then we may apply \ref{Sah7.2} and get $s_0,s_1\in\Sigma^*(
t)$ such that $\nor[s_\ell]=\nor[t]-1$, $\pos(s_\ell)\subseteq\{\nu\in\pos(
t):r^\ell_\nu>0\}$ and 
\[c_\ell\stackrel{\rm def}{=}F^*_{s_\ell}(r^\ell_\nu:\nu\in\pos(s_\ell))
\geq a_\ell\cdot (1-2^{-2^{k+3}}).\]
Then 
\[c_0+c_1\geq (a_0+a_1)\cdot (1-2^{-2^{k+3}})\geq F^*_t(r_\nu:\nu\in\pos(t
))\cdot (1-2^{-2^k}),\]
and we are done.

So suppose now that $a_\ell<2^{-2^{k+3}}$. Then 
\[a_{1-\ell}\geq 2^{|g_t|-N_k}\cdot\sum\{r_\nu:\nu\in\pos(t)\}-2^{-2^{k+3}}
\geq 2^{-2^k}-2^{-2^{k+3}}\geq 2^{-2^{k+3}},\]
and using \ref{Sah7.2} we find $s_{1-\ell}\in\Sigma^*(t)$ such that
$\nor[s_{1-\ell}]=\nor[t]-1$, $\pos(s_{1-\ell})\subseteq\{\nu\in\pos(t):
r^{1-\ell}_\nu>0\}$, and 
\[\begin{array}{l}
c_{1-\ell}\stackrel{\rm def}{=}F^*_{s_{1-\ell}}(r^{1-\ell}_\nu:\nu\in\pos(
s_{1-\ell}))\geq a_{1-\ell}\cdot (1-2^{-2^{k+3}})\geq{}\\
(F^*_t(r_\nu:\nu\in\pos(t))-2^{-2^{k+3}})\cdot(1-2^{-2^{k+3}})\geq{}\\
F^*_t(r_\nu:\nu\in\pos(t))\cdot (1-2^{-2^k})+2^{-2^k}(2^{-2^k}-2^{-2^{k+3}})
-2^{-2^{k+3}}+2^{-2^{k+4}}={}\\
F^*_t(r_\nu:\nu\in\pos(t))\cdot (1-2^{-2^k})+2^{-2^{k+1}}-2^{-9\cdot 2^k}-
2^{-2^{k+3}}+2^{-2^{k+4}}\geq{}\\
F^*_t(r_\nu:\nu\in\pos(t))\cdot(1-2^{-2^k}). 
\end{array}\]
\end{proof}

The following lemma and the proposition are, as a matter of fact, included
in \ref{getmore}, \ref{crucial}. However, we decided that \ref{prenullpres}
and \ref{nullpres} could be a good warm-up, and also we will use their
proofs later.

\begin{lemma}
\label{prenullpres}
Assume that:
\begin{enumerate}
\item[(i)] $t\in K^*$, $\nor[t]>1$, $k=k_t$, $\gamma\in [0,1]$,
\item[(ii)] $\langle r_\nu:\in\pos(t)\rangle\subseteq [0,1]$, $a=
F^*_t(r_\nu:\nu\in\pos(t))$, $\gamma\cdot a\geq 2^{-6\cdot 2^k}$,
\item[(iii)] $Y$ is a finite set,
\item[(iv)] for $\nu\in \pos(t)$, $u_\nu$ is a function from $Y$ to $[0,1]$
such that  
\[\gamma\cdot r_\nu\cdot |Y|\leq \sum\{u_\nu(y):y\in Y\},\]
\item[(v)] for $y\in Y$ we let
\[\begin{array}{ll}
u(y)=\sup\{b:&\mbox{there is }s\in\Sigma^*(t)\mbox{ such that }\nor[s]\geq
\nor[t]-1\mbox{ and}\\
&b\leq F^*_s(u_\nu(y):\nu\in\pos(s))\}.
  \end{array}\]
\end{enumerate}
Then
\[\gamma\cdot a\cdot (1-2^{-2^k})\leq\frac{\sum\{u(y):y\in Y\}}{|Y|}.\]
\end{lemma}

\begin{proof}
Let $k=k_t$, $N=N_{k_t}$, $g=g_t$.

First note that
\[\begin{array}{l}
a=F^*_t(r_\nu:\nu\in\pos(t))\leq 2^{|g|-N}\cdot\sum\{r_\nu:\nu\in\pos(t)\}
\leq{}\\
2^{|g|-N}\cdot\frac{1}{\gamma}\cdot\frac{1}{|Y|}\cdot\sum\limits_{\nu
\in \pos(t)}\sum\limits_{y\in Y} u_\nu(y)=
\frac{1}{\gamma}\cdot\frac{1}{|Y|}\cdot\sum\limits_{y\in Y}\left(2^{|g|-N}
\cdot\sum\limits_{\nu\in\pos(t)} u_\nu(y)\right).
  \end{array}\]
Let $C\stackrel{\rm def}{=}\{y\in Y:2^{|g|-N}\cdot\sum\limits_{\nu\in\pos(
t)}  u_\nu(y)\geq 2^{-2^{k+3}}\}$. For each $y\in C$ we may use \ref{Sah7.2}
to pick $s_y\in\Sigma^*(t)$ such that $\nor[s_y]\geq \nor[t]-1$ and 
\[F^*_{s_y}(u_\nu(y):\nu\in\pos(s_y))\geq 2^{|g|-N}\cdot\sum\limits_{\nu\in
\pos(t)} u_\nu(y)\cdot (1-2^{-2^{k+3}}).\]
Hence,
\[\begin{array}{l}
a\leq\frac{1}{\gamma}\cdot \frac{|Y\setminus C|}{|Y|}\cdot 2^{-2^{k+3}}+
\frac{1}{\gamma}\cdot\frac{1}{|Y|}\cdot\sum\limits_{y\in C}
\frac{F^*_s(u_\nu(y):\nu\in\pos(s_y))}{1-2^{-2^{k+3}}}\leq{}\\
\frac{1}{\gamma}\cdot 2^{-2^{k+3}}+\frac{1}{\gamma}\cdot\frac{1}{|Y|}\cdot
\frac{1}{1-2^{-2^{k+3}}}\cdot\sum\limits_{y\in C}u(y).
  \end{array}\]
Consequently,
\[(\gamma a-2^{-2^{k+3}})(1-2^{-2^{k+3}})\leq\frac{\sum\limits_{y\in C}
u(y)}{|Y|}\leq \frac{\sum\limits_{y\in Y} u(y)}{|Y|},\]
and hence
\[\gamma a(1-2^{-2^k})\leq\frac{\sum\limits_{y\in Y} u(y)}{|Y|}.\]
\end{proof}

\begin{proposition}
\label{nullpres}
The forcing notion $\starQ$ preserves outer (Lebesgue) measure.
\end{proposition}

\begin{proof}
Assume that $A\subseteq\prod\limits_{i<\omega} N_i$ is a set of outer
(Lebesgue) measure $1$. We have to show that, in $\bV^\starQ$, it is still
an outer measure one set. 

Let $\dot{T}$ be a $\starQ$--name for a tree such that $\dot{T}\subseteq
\bigcup\limits_{i\in\omega}\prod\limits_{j<i}N_i$ and the Lebesgue measure
$\mu([\dot{T}])$ of the set $[\dot{T}]$ of $\omega$--branches through
$\dot{T}$ is positive, and suppose that some condition $p$ forces
``$[\dot{T}]\cap A=\emptyset$''. Take a condition $q\geq p$ such that 
\begin{enumerate}
\item[$(\alpha)$] $q$ is special (remember \ref{nordense}) and $\lh(\mrot(
q))=k_0>5$, and $\nor[t^q_\eta]>2$ for all $\eta\in T^q$, and $\mu^{\bF^*}(
q)>\frac{1}{2}$, 
\item[$(\beta)$]  for some $\rho\in\prod\limits_{j<n} N_j$, $n<\omega$, the
condition $q$ forces that $\mu([(\dot{T})^{[\rho]}])\cdot \prod\limits_{j<n}
N_j\geq \frac{7}{8}$, 
\item[$(\gamma)$] for some $k_0<k_1<k_2<\ldots$, letting $F_i=T^q\cap
\prod\limits_{m<k_i}\bH^*(m)$, we have that for each $\nu\in F_i$, the
condition $q^{[\nu]}$ decides the value of $\dot{T}\cap\prod\limits_{j<n+i}
N_j$ (remember \ref{contreading}).  
\end{enumerate}
Fix $i<\omega$ for a moment, and let $Y_i=\{y\in\prod\limits_{j<n+i}
N_j:\rho\vartriangleleft y\}$.

For $\nu\in T[q,F_i]$ and $y\in Y_i$ we let 
\[\begin{array}{ll}
u_\nu(y)=\sup\{\mu^{\bF^*}(q'):&q'\mbox{ is a condition stronger than $q$
and such that}\\
&\mrot(q')=\nu\ \mbox{ and }\ (\forall\eta\in T^{q'})(\nor[t^{q'}_\eta]\geq
\nor[t^q_\eta]-1)\\
&\mbox{and }\ q'\forces y\in\dot{T}\}. 
  \end{array}\]

\begin{claim}
\label{cl2}
If $\eta\in T[q,F_i]$, $k_0\leq\lh(\eta)=k\leq k_i$, then
\[\frac{7}{8}\cdot\prod\limits_{\ell=k}^{k_i-1}(1-2^{-2^\ell})\cdot|Y_i|\cdot 
\mu^{\bF^*}_q(\eta)\leq \sum\big\{u_\eta(y):y\in Y_i\big\}.\]  
[If $k=k_i$, then we stipulate $\prod\limits_{\ell=k}^{k_i-1}(1-2^{-2^\ell})
=1$.] 
\end{claim}

\begin{proof}[Proof of the claim]
We show it by downward induction on $\eta\in T[q,F_i]$. If $k=\lh(\eta)=
k_i$, then $q^{[\eta]}$ decides $\dot{T}\cap Y_i$, and if $q^{[\eta]}$
forces that $y\in\dot{T}\cap Y_i$, then $u_\eta(y)\geq\mu^{\bF^*}_q(\eta)$. 
Hence, by $(\beta)$, we have $\frac{7}{8}\cdot |Y_i|\cdot\mu^{\bF^*}_q(\eta)
\leq\sum\{u_\eta(y):y\in Y_i\}$. 

Let us assume now that $k=\lh(\eta)=k_i-1$. Apply \ref{prenullpres} to
$t^q_\eta$, $\gamma=\frac{7}{8}$, $u_\nu$, $r_\nu=\mu^{\bF^*}_q(\nu)$ (for
$\nu\in\pos(t^q_\eta)$), and $Y_i$. Note that, as $q$ is special,
$\mu^{\bF^*}_q(\eta)\geq 2^{-2^{k+1}}$, so $\gamma\cdot F^*_{t^q_\eta}(
r_\nu:\nu\in\pos(t^q_\eta))=\frac{7}{8}\mu^{\bF^*}_q(\eta)>2^{-6\cdot
2^k}$. Also note that 
\begin{enumerate}
\item[$(*)$] $u(y)$ defined as in \ref{prenullpres}(v) is $u_\eta(y)$. 
\end{enumerate}
[Why? First suppose that $u(y)<u_\eta(y)$. By the definition of $u_\eta$ we
may find $q'\geq q$ such that $\mrot(q')=\eta$, $\nor[t^{q'}_\nu]\geq
\nor[t^q_\nu]-1$ for $\nu\in T^{q'}$, and $q'\forces y\in\dot{T}$, and
$\mu^{\bF^*}(q')>u(y)$. Note that $\mu^{\bF^*}_{q'}(\nu)\leq u_\nu(y)$ for
all $\nu\in\pos(t^{q'}_\eta)$, and thus
\[u(y)<\mu^{\bF^*}(q')=F^*_{t^{q'}_\eta}(\mu^{\bF^*}_{q'}(\nu):\nu\in\pos(
t^{q'}_\eta))\leq F^*_{t^{q'}_\eta}(u_\nu(y):\nu\in\pos(t^{q'}_\eta)).\]
By the definition of $u(y)$, the last expression is $\leq u(y)$, a
contradiction. Now suppose $u(y)>u_\eta(y)$. Take $s\in\Sigma^*(t^q_\eta)$
such that $\nor[s]\geq\nor[t^q_\eta]-1$ and $F^*_s(u_\nu(y):\nu\in\pos(s))
>u_\eta(y)$; clearly we may request that $u_\nu(y)>0$ for $\nu\in
\pos(s)$. Let $z_\nu<u_\nu(y)$ (for $\nu\in\pos(s)$) be positive numbers
such that if $z_\nu\leq r_\nu\leq u_\nu(y)$ for $\nu\in\pos(s)$, then
$F^*_s(r_\nu:\nu\in\pos(s))>u_\eta(y)$ (compare \ref{easyprop}). Pick
conditions $q'_\nu$ such that $\mu^{\bF^*}(q'_\nu)>z_\nu$, $q'_\nu$ as in
definition of $u_\nu(y)$, and let $q'$ be such that $\mrot(q')=\eta$,
$t^{q'}_\eta=s$, and $(q')^{[\nu]}=q'_\nu$ for $\nu\in\pos(s)$. Then
$\mu^{\bF^*}(q')>u_\eta(y)$ giving an easy contradiction.]\\
Thus we get 
\[\frac{7}{8}\cdot \mu^{\bF^*}_q(\eta)\cdot (1-2^{-2^{k_i-1}})\cdot |Y_i|\leq
\sum\{u_\eta(y):y\in Y_i\},\]
as required.

Now suppose $k_0\leq k=\lh(\eta)<k_i-1$, and we have proved the assertion of
the claim for all $\nu\in\pos(t^q_\eta)$. We again apply \ref{prenullpres},
this time to $\gamma=\frac{7}{8}\cdot\prod\limits_{\ell=k+1}^{k_i-1}(1-
2^{-2^\ell})$, and $t^q_\eta$, $u_\nu$, $r_\nu=\mu^{\bF^*}_q(\nu)$ (for
$\nu\in\pos(t^q_\eta)$) and $Y_i$. We note that 
\[\begin{array}{r}
\frac{7}{8}\cdot\prod\limits_{\ell=k+1}^{k_i-1}(1-2^{-2^\ell})\cdot
F^*_{t^q_\eta}(r_\nu:\nu\in\pos(t^q_\eta))=\frac{7}{8}\cdot
\prod\limits_{\ell=k+1}^{k_i-1}(1-2^{-2^\ell})\cdot\mu^{\bF^*}_q(\eta)\geq
{}\\
\frac{7}{8}\cdot (1-2^{1-2^{k+1}})\cdot 2^{-2^{k+1}}\geq 2^{-6\cdot 2^k},
  \end{array}\]
so the assumptions of \ref{prenullpres} are satisfied. Therefore we may
conclude that
\[\frac{7}{8}\cdot\prod\limits_{\ell=k+1}^{k_i-1}(1-2^{-2^\ell})\cdot
\mu^{\bF^*}_q(\eta)\cdot (1-2^{-2^k})\cdot |Y_i|\leq \sum\{u_\eta(y): y\in
Y_i\},\]
as needed.
\end{proof}
Applying \ref{cl2} to $\eta=\mrot(q)$ we get 
\[\frac{7}{8}\cdot\prod\limits_{\ell=k_0}^{k_i-1}(1-2^{-2^\ell})\cdot
\mu^{\bF^*}(q)\leq \frac{\sum\big\{u_{\mrot(q)}(y):y\in Y_i\big\}}{|Y_i|},\] 
and hence $\frac{3}{4}\mu^{\bF^*}(q)\cdot |Y_i|\leq \sum\{u_{\mrot(q)}(y):
y\in Y_i\}$. Then necessarily 
\[\frac{1}{4}|Y_i|\leq |\{y\in Y_i:u_{\mrot(p)}(y)\geq
\frac{1}{4}\mu^{\bF^*} (q)\}|\]
(remember $\mu^{\bF^*}(q)>\frac{1}{2}$). Let $Z_i= \{y\in Y_i:u_{\mrot(p)}
(y)\geq \frac{1}{4}\mu^{\bF^*}(q)\}$. 
\medskip

Look at the set $\{x\in\prod\limits_{i<\omega} N_i:(\exists^\infty i<\omega)
(x\restriction (n+i)\in Z_i)\}$ -- it is a Borel set of positive (Lebesgue)
measure, and therefore we may pick $x\in A$ such that $(\exists^\infty
i<\omega)(x\restriction (n+i)\in Z_i)$. For each $i<\omega$ such that
$x\restriction (n+i)\in Z_i$ choose a condition $q_i\in\starQ$ such that 
\begin{itemize}
\item $q_i\geq q$, $\mrot(q_i)=\mrot(q)$, $\mu^{\bF^*}(q_i)>\frac{1}{8}
\mu^{\bF^*}(q)$, and
\item $(\forall\eta\in T^{q_i})(\nor[t^{q_i}_\eta]\geq \nor[t^q_\eta]-1)$,
and 
\item $q_i\forces x\restriction (n+i)\in \dot{T}$.
\end{itemize}
By K\"onig Lemma (remember $(K^*,\Sigma^*)$ is strongly finitary) we find an
infinite set $I\subseteq\omega$ such that for each $i<j_0<j_1$ from $I$ we
have  
\[T^{q_{j_0}}\cap\!\prod_{k<k_i}\!\bH(k)=T^{q_{j_1}}\cap\!\prod_{k<k_i}\!
\bH(k)\quad\mbox{and}\quad (\forall\eta\in T^{q_{j_0}})(\lh(\eta)<k_i\
\Rightarrow\ t^{q_{j_0}}_\eta=t^{q_{j_1}}_\eta).\]
Let $q^*=\langle s_\eta:\eta\in S\rangle$ be such that $\mrot(S)=\mrot(q)$,
\[S=\bigcup\limits_{i\in I}\{\eta\in T^{q_j}:j\in I\ \&\ i<j\ \&\ \lh(\eta)
< k_i\},\]
and if $\eta\in S$, then $\suc_S(\eta)=\pos(s_\eta)$ and
$s_\eta=t^{q_i}_\eta$  for sufficiently large $i\in I$. It should be clear
that $q^*\in\starQ$ is a condition stronger than $q$, and it forces that
$x\in [\dot{T}]\cap A$, a contradiction.  
\end{proof}

\begin{lemma}
\label{getmore}
Assume that:
\begin{enumerate}
\item[(i)] $t\in K^*$, $\nor[t]>1$, $k=k_t>1$, $\gamma\in [0,1]$,
\item[(ii)] $\langle r_\nu:\in\pos(t)\rangle\subseteq [0,1]$, $a=
F^*_t(r_\nu:\nu\in\pos(t))$, $\gamma\cdot a\geq 2^{-6\cdot 2^k}$,
\item[(iii)] $Y^*$ is a finite set, $Y= Y^*\times N_k$,
\item[(iv)] for $\nu\in \pos(t)$, $u_\nu$ is a function from $Y$ to $[0,1]$
such that  
\[\gamma\cdot r_\nu\cdot |Y|\leq \sum\{u_\nu(y):y\in Y\},\]
\item[(v)] for $y=(y_0,y_1)\in Y^*\times N_k$ and $\ell<2$ we let
\[\begin{array}{ll}
u(y,\ell)=\sup\{b:&\mbox{there is }s\in\Sigma^*(t)\mbox{ such that }\nor[s] 
\geq\nor[t]-1\mbox{ and}\\
&(\forall\nu\in\pos(s))(\nu(k)(y_1)=\ell)\mbox{ and }b\leq F^*_s(u_\nu(y):
\nu\in\pos(s))\}.
  \end{array}\]
\end{enumerate}
Then
\[\gamma\cdot a\cdot (1-2^{-2^k})\leq\frac{1}{2\cdot|Y|}\sum\{u(y,\ell):y\in
Y\ \&\ \ell<2\}.\]
\end{lemma}

\begin{proof}
Let $k=k_t$, $N=N_k$, $g=g_t$. Note that 
\[\begin{array}{l}
a=F^*_t(r_\nu:\nu\in\pos(t))\leq 2^{|g|-N}\cdot\sum\{r_\nu:\nu\in\pos(t)\}
\leq{}\\
2^{|g|-N}\cdot\frac{1}{\gamma}\cdot\frac{1}{|Y|}\cdot\sum\limits_{\nu
\in \pos(t)}\sum\limits_{\ell<2}\Big(\sum \{u_\nu(y_0,y_1): (y_0,y_1)\in Y\
\&\ \nu(k)(y_1)=\ell\}\Big)=\\
\frac{1}{\gamma}\cdot\frac{1}{2\cdot |Y|}\cdot\sum\limits_{(y_0,y_1,\ell)\in
Y\times 2}\Big(2^{|g|-N+1}\cdot\sum \{u_\nu(y_0,y_1):\nu\in\pos(t)\ \&\
\nu(k)(y_1)=\ell\}\Big). 
  \end{array}\]
Let $C$ consist of all triples $(y_0,y_1,\ell)\in Y^*\times N\times 2$ such
that $y_1\notin\dom(g)$ and 
\[2^{|g|+1-N}\cdot\sum\{u_\nu(y_0,y_1):\nu\in\pos(t)\ \&\ \nu(k)(y_1)=\ell\}
\geq 2^{-2^{k+3}},\]
and fix $(y_0,y_1,\ell)\in C$ for a moment. Let $g':\dom(g)\cup\{y_1\}
\longrightarrow 2$ be such that $g\subseteq g'$ and $g'(y_1)=\ell$. Apply
\ref{Sah7.2} (to $t,g'$ and $u_\nu(y_0,y_1)$ for $\nu\in\pos(t)$,
$g'\subseteq\nu(k)$) to pick $s=s_{y_0,y_1,\ell}\in\Sigma^*(t)$ such that
$\nor[s]\geq \nor[t]-1$, $g'\subseteq g_s$ and   
\[\frac{F^*_s(u_\nu(y_0,y_1):\nu\in\pos(s))}{1-2^{-2^{k+3}}}\geq 2^{|g|+1-N} 
\cdot\sum\{u_\nu(y_0,y_1):\nu\in\pos(t)\ \&\ g'\subseteq\nu(k)\}.\]
Next note that $\frac{|g|}{N}<2^{-2^{k+3}}$, so 
\[\begin{array}{r}
\frac{1}{\gamma\cdot |Y|}\cdot\sum\limits_{(y_0,y_1,\ell)\in
Y\times 2\setminus C}\Big(2^{|g|-N}\cdot\sum \{u_\nu(y_0,y_1):\nu\in\pos(t)\
\&\ \nu(k)(y_1)=\ell\}\Big)\leq{}\\
\frac{|g|}{\gamma\cdot N}+\frac{1}{\gamma}\cdot 2^{-2^{k+3}}\leq
\frac{1}{\gamma}\cdot 2^{1-2^{k+3}}.
  \end{array}\]
Therefore,
\[\begin{array}{l}
a\leq\frac{1}{\gamma}\cdot 2^{1-2^{k+3}}+\frac{1}{\gamma}\cdot\frac{1}{2
\cdot |Y|}\cdot\sum\limits_{(y_0,y_1,\ell)\in C}\frac{F^*_{s_{y_0,y_1,\ell}}
(u_\nu(y_0,y_1):\nu\in\pos(s_{y_0,y_1,\ell}))}{1-2^{-2^{k+3}}}\leq{}\\
\frac{1}{\gamma}\cdot 2^{1-2^{k+3}}+\frac{1}{\gamma}\cdot\frac{1}{2\cdot
|Y|}\cdot\frac{1}{1-2^{-2^{k+3}}}\cdot\sum\limits_{(y,\ell)\in C}u(y,\ell).
  \end{array}\]
Hence,
\[(\gamma a-2^{1-2^{k+3}})(1-2^{-2^{k+3}})\leq\frac{1}{|Y\times 2|}
\sum\limits_{(y,\ell)\in Y\times 2}u(y,\ell),\]
and therefore, as $\gamma a\geq 2^{-6\cdot 2^k}$ and $k>1$, 
\[\gamma a(1-2^{-2^k})\leq\frac{1}{|Y\times 2|}\sum\limits_{(y,\ell)\in
Y\times 2} u(y,\ell).\]
\end{proof}

Let $\dot{W}$ be the canonical $\starQ$--name for the generic real (so
$\dot{W}$ is a name for a function in $\prod\limits_{i<\omega}\bH^*(i)$ such
that $p\forces\mrot(p)\subseteq\dot{W}$). Also, let $\dot{h}$ be a name for
a function from $\prod\limits_{i<\omega} N_i$ to $\can$ such that
$\dot{h}(x)(i)=\dot{W}(i)\big(x(i)\big)$. Clearly, $\dot{h}$ is (a name for)
a continuous function.

Now comes the main property of the forcing notion $\starQ$.

\begin{proposition}
\label{crucial}
Suppose that $A\subseteq\prod\limits_{i<\omega}N_i\times \can$ is a set of
outer (Lebesgue) measure 1. Then, in $\bV^{\starQ}$, the set 
\[\{x\in\prod\limits_{i<\omega} N_i:(x,\dot{h}(x))\in A\}\]
has outer measure 1.
\end{proposition}

\begin{proof}
Assume toward contradiction that $\dot{T}$ is a $\starQ$--name for a tree
included in $\bigcup\limits_{k<\omega}\prod\limits_{i<k} N_i$, and
$p\in\starQ$ is a condition such that 
\[p\forces_{\starQ}\mbox{`` }\mu([\dot{T}])>0\mbox{ and }(\forall x\in
[\dot{T}])((x,\dot{h}(x))\notin A)\mbox{ ''}.\] 
(Here, $\mu$ stands for the product measure on $\prod\limits_{i<\omega}
N_i$.) Passing to a stronger condition and shrinking the tree $\dot{T}$ (if 
necessary) we may assume that  
\begin{enumerate}
\item[$(\alpha)$] $p$ is special and $\lh(\mrot(p))=k_0>5$, and
$\nor[t^p_\eta]>2$ for all $\eta\in T^p$, and $\mu^{\bF^*}(p)>\frac{1}{2}$, 
\item[$(\beta)$]  for some $\rho\in\prod\limits_{j<n} N_j$, $n<k_0$, the
condition $p$ forces that $\mu([(\dot{T})^{[\rho]}])\cdot \prod\limits_{j<n}
N_j\geq \frac{7}{8}$, 
\item[$(\gamma)$] for some $k_0<k_1<k_2<\ldots$, letting $F_i=T^p\cap
\prod\limits_{m<k_i}\bH^*(m)$, we have that for each $\nu\in F_{i+1}$, the
condition $p^{[\nu]}$ decides the value of $\dot{T}\cap\prod\limits_{j<k_i}
N_j$.  
\end{enumerate}
Fix $i<\omega$ for a moment, and let $Y_i^{**}=\{y\in\prod\limits_{j<k_i}
N_j:\rho\vartriangleleft y\}$.

Let $\nu_0\in F_i$, and for $\nu\in T[p^{[\nu_0]},F_{i+1}]$ and $y\in
Y^{**}_i$ let
\[\begin{array}{ll}
u_\nu(y)=\sup\{\mu^{\bF^*}(p'):&p'\mbox{ is a condition stronger than $p$
and such that}\\
&\mrot(p')=\nu\ \mbox{ and }\ (\forall\eta\in T^{p'})(\nor[t^{p'}_\eta]\geq 
\nor[t^p_\eta]-1),\\
&\mbox{and }\ p'\forces y\in\dot{T}\}. 
  \end{array}\]
So we are at the situation from the proof of \ref{nullpres} (with $q$ there 
replaced by $p$), and we may use \ref{cl2} to conclude that
\begin{enumerate}
\item[$(\circledast)$] \quad $\frac{7}{8}\cdot\prod\limits_{\ell=k_i}^{k_{i
+1}-1}(1-2^{-2^\ell})\cdot |Y^{**}_i|\cdot\mu^{\bF^*}_p(\nu_0)\leq\sum
\{u_{\nu_0}(y): y\in Y^{**}_i\}$. 
\end{enumerate}
Now, for each $\nu\in T[p,F_i]$ we define $u^*_\nu:Y^{**}_i\times
2^{\textstyle [\lh(\nu),k_i)}\longrightarrow [0,1]$ by
\[\begin{array}{ll}
u^*_\nu(y,\sigma)=&\\
\sup\{\mu^{\bF^*}(p'):&p'\mbox{ is a condition stronger than $p$ and such
that}\\ 
&\mrot(p')=\nu\ \mbox{ and }\ (\forall\eta\in T^{p'})(\nor[t^{p'}_\eta]\geq 
\nor[t^p_\eta]-1),\\
&\mbox{and }\ p'\forces\mbox{`` } y\in\dot{T}\ \&\ \Big(\forall j\in
[\lh(\nu),k_i)\Big)\Big(\dot{W}(j)\big(y(j)\big)=\sigma(j)\Big)\mbox{ ''}\}.  
  \end{array}\]
(If $\nu\in F_i$, so $\lh(\nu)=k_i$, then $2^{\textstyle[\lh(\nu),k_i)}=
\{\emptyset\}$ and $u^*_\nu(y,\emptyset)=u_\nu(y)$.)

\begin{claim}
\label{cl3}
If $\eta\in T[p,F_i]$, $k_0\leq\lh(\eta)=k\leq k_i$, then 
\[\frac{7}{8}\cdot\prod\limits_{\ell=k}^{k_{i+1}-1}(1-2^{-2^\ell})\cdot
|Y^{**}_i|\cdot 2^{k_i-k}\cdot\mu^{\bF^*}_p(\eta)\leq\sum\{u_\eta^*(y,
 \sigma): (y,\sigma)\in X^i_\eta\},\]
where $X^i_\eta=Y^{**}_i\times 2^{\textstyle[\lh(\eta),k_i)}$. 
\end{claim}

\begin{proof}[Proof of the claim]
The proof, by downward induction on $\eta$, is similar to that of \ref{cl2},
but this time we use \ref{getmore}.

First note that if $k=k_i$, then our assertion is exactly what is stated in
$(\circledast)$. So suppose that $\eta\in T[p,F_i]$, $\lh(\eta)<k_i$, and
that we have proved our claim for all $\nu\in\pos(t^p_\eta)$. We are going
to apply \ref{getmore} to $t=t^p_\eta$, $\gamma=\frac{7}{8}\cdot
\prod\limits_{\ell=k+1}^{k_{i+1}-1}(1-2^{-2^\ell})$, $Y^*=\{y\restriction
(k_i\setminus\{k\}):y\in Y^{**}_i\}\times 2^{\textstyle [k+1,k_i)}$ (and
$Y=Y^*\times N_k$ being interpreted as $Y^{**}_i\times 2^{\textstyle [k+1,
k_i)}$), and $r_\nu=\mu^{\bF^*}_p(\nu)$, and $u_\nu(y,\sigma)=u^*_\nu(y,
\sigma)$ (for $\nu\in\pos(t^p_\eta)$, $(y,\sigma)\in X^i_\nu$), so we have
to check the assumptions there. Note that (as $p$ is special)
\[\gamma\cdot F^*_t(r_\nu:\nu\in\pos(t))=\gamma\cdot \mu^{\bF^*}_p(\eta)
\geq\frac{7}{8}\cdot\prod_{\ell=k+1}^{k_{i+1}-1}(1-2^{-2^\ell})\cdot
2^{-2^{k+1}}>2^{-6\cdot 2^k}\]
(so the demand in \ref{getmore}(ii) is satisfied). Also, by the inductive
hypothesis, 
\[\gamma\cdot |Y^*\times N_k|\cdot r_\nu\leq \sum\{u^*_\nu(y,\sigma):(y,
\sigma)\in X^i_\nu\}\]
(so \ref{getmore}(iv) holds). Finally note that if $(y,\sigma)\in Y^{**}_i
\times 2^{\textstyle [k+1,k_i)}$, $\ell<2$, and $\sigma':[k,k_i)
\longrightarrow 2$ is such that $\sigma'(k)=\ell$, $\sigma\restriction [k+1,
k_i)=\sigma$, then $u(y,\sigma,\ell)$ defined by \ref{getmore}(v) is
$u^*_\eta(y,\sigma')$.  

So, by \ref{getmore}, we may conclude that
\[\begin{array}{r}
\frac{7}{8}\cdot\prod\limits_{\ell=k+1}^{k_{i+1}-1}(1-2^{-2^\ell})\cdot
\mu^{\bF^*}_p(\eta)\cdot (1-2^{-2^k})\cdot 2\cdot |Y^{**}_i|\cdot
2^{k_i-k-1}\leq{}\quad\\
\sum\{u_\eta^*(y,\sigma'): (y,\sigma')\in X^i_\eta\},
  \end{array}\] 
as needed.
\end{proof}
In particular, it follows from \ref{cl3} that
\[\frac{7}{8}\cdot\prod\limits_{\ell=k_0}^{k_{i+1}-1}(1-2^{-2^\ell})\cdot
\mu^{\bF^*}(p)\leq\frac{\sum\{u_{\mrot(p)}^*(y,\sigma): (y,\sigma)\in
Y^{**}_i\times 2^{\textstyle [k_0,k_i)}\}}{|Y^{**}_i|\cdot 2^{k_i-k_0}},\] 
and hence 
\[\frac{3}{4}\cdot\mu^{\bF^*}(p)\leq\frac{\sum\{u_{\mrot(p)}^*(y,\sigma):
(y,\sigma)\in Y^{**}_i\times 2^{\textstyle [k_0,k_i)}\}}{|Y^{**}_i|\cdot
2^{k_i-k_0}}.\]  
Let $\varphi:\prod\limits_{j<k_0}N_j\longrightarrow 2^{\textstyle k_0}$ be
such that $\varphi(y)(j)=(\mrot(p)(j))(y(j))$. Now we define:
\[\begin{array}{lcl}
Z_i&=&\{(y,\sigma)\in Y^{**}_i\times 2^{\textstyle [k_0,k_i)}:
u^*_{\mrot(p)}(y,\sigma)\geq \frac{1}{4}\mu^{\bF^*}(p)\},\quad\mbox{ and}\\
Z_i^+&=&\{(y,\sigma)\in Y^{**}_i\times 2^{\textstyle k_i}:\varphi(y
\restriction k_0)=\sigma\restriction k_0\,\ \&\ \, (y,\sigma\restriction
[k_0, k_i))\in Z_i\}. 
\end{array}\]
Note that $|Z_i|\geq \frac{1}{4}|Y^{**}_i\times 2^{\textstyle [k_0,k_i)}|$, 
and therefore
\[\frac{|Z^+_i|}{2^{k_i}\cdot\prod\limits_{j<k_i}N_j}\geq \frac{1}{2^{k_0
+2}\cdot\prod\limits_{j<k_0} N_j}.\]
Now we may finish like in \ref{nullpres}: the set
\[\{(x_0,x_1)\in\prod_{j<\omega}N_j\times\can:(\exists^\infty i<\omega)(
(x_0\rest k_i,x_1\rest k_i)\in Z^+_i)\}\]
is a Borel set of positive (Lebesgue) measure, so we may choose $(x_0,x_1)
\in A$ such that for infinitely many $i<\omega$ we have $(x_0\rest k_i,x_1
\rest k_i)\in Z^+_i$. For each such $i$ pick a condition $q_i\geq p$ such
that 
\begin{itemize}
\item $\mrot(q_i)=\mrot(p)$, $\mu^{\bF^*}(q_i)>\frac{1}{8}\mu^{\bF^*}(p)$,
and 
\item $(\forall\eta\in T^{q_i})(\nor[t^{q_i}_\eta]\geq\nor[t^p_\eta]-1)$,
and  
\item $q_i\forces\mbox{`` } x_0\restriction k_i\in \dot{T}\mbox{ and }
\Big(\forall j\in [k_0,k_i)\Big)\Big(\dot{W}(j)\big(x_0(j)\big)=x_1(j)\Big)
\mbox{ ''}$.
\end{itemize}
By K\"onig Lemma, we may find a condition $q\in\starQ$ stronger than $p$,
and an infinite set $I\subseteq\omega$ such that 
\begin{enumerate}
\item[$(\otimes)$] if $i<j$ are from $I$, then $i+1<j$ and
\[T^{q_j}\cap\prod_{k<k_{i+1}} N_k =T^q\cap\prod_{k<k_{i+1}} N_k\quad
\mbox{and}\quad (\forall\eta\in T^{q_j})(\lh(\eta)<k_{i+1}\ \Rightarrow\
t^{q_j}_\eta=t^q_\eta).\]
\end{enumerate}
Then clearly $q\forces\mbox{`` }x_0\in\dot{T}\ \&\ \dot{h}(x_0)=x_1\mbox{
''}$, a contradiction.
\end{proof}

\section{The first model: sup-measurability}
To prove the first of our main results, let us start with a reduction of
the sup-measurability problem.
\begin{lemma}
\label{boxtimes}
The following conditions are equivalent:
\begin{enumerate}
\item[$(\boxtimes)_{\rm sup}^1$] Every sup-measurable function
$f:\mbR\times\mbR\longrightarrow\mbR$ is Lebesgue measurable.
\item[$(\boxtimes)_{\rm sup}^2$] For every non-measurable set $A\subseteq
\mbR\times\mbR$ there exists a Borel function $f:\mbR\longrightarrow\mbR$ such
that the set $\{x\in\mbR:(x,f(x))\in A\}$ is not measurable.
\item[$(\boxtimes)_{\rm sup}^3$] For every non-measurable set $A\subseteq
\can\times\can$ there is a Borel function $f:\can\longrightarrow\can$ such
that the set $\{x\in\can:(x,f(x))\in A\}$ is not measurable.
\item[$(\boxtimes)_{\rm sup}^4$] For every set $A\subseteq\prod\limits_{k<
\omega}N_k\times\can$ of outer measure one and inner measure zero, there is
a Borel function $h:\prod\limits_{k<\omega}N_k\longrightarrow\can$ such that
the set 
\[\{x\in\prod\limits_{k<\omega}N_k:(x,h(x))\in A\}\]
is not measurable.\\ 
(Here, the sequence $\langle N_k:k<\omega\rangle$ is the one defined at the
beginning of the second section.)
\end{enumerate}
\end{lemma}

\begin{proof}
The equivalences $(\boxtimes)_{\rm sup}^1\Leftrightarrow (\boxtimes)_{\rm
sup}^2\Leftrightarrow(\boxtimes)_{\rm sup}^3$ are well known (see Balcerzak
\cite[Proposition 1.5]{Ba91b}; also compare with the proof of Ciesielski and 
Shelah \cite[Corollary 3]{CiSh:695}). 
\medskip

\noindent $(\boxtimes)_{\rm sup}^4\Rightarrow(\boxtimes)_{\rm sup}^3$:\qquad
Assume $(\boxtimes)_{\rm sup}^4$, and suppose that $A\subseteq\can\times
\can$ is a non-measurable set. Then we may find a closed set $C\subseteq
\can\times\can$ such that
\begin{itemize}
\item for each $x\in\can$, the set $\{y\in\can:(x,y)\in C\}$ is either empty
or is a perfect set of positive Lebesgue measure,
\item for every Borel set $D\subseteq C$ of positive measure, both $A\cap
D\neq\emptyset$ and $D\setminus A\neq\emptyset$ (that is, both $A\cap C$ and
$C\setminus A$ are of full outer measure in $C$).
\end{itemize}
Pick a Borel isomorphism $\psi=(\psi_0,\psi_1):C\longrightarrow\prod\limits_{k<
\omega}N_k\times\can$ such that
\begin{itemize}
\item if $(x,y),(x',y')\in C$, then $\psi_0(x,y)=\psi_0(x',y')\
\Leftrightarrow\ x=x'$, 
\item if $B\subseteq C$ is Borel, then $B$ has measure 0 if and only if its
image $\psi[B]$ has measure zero. 
\end{itemize}
Now note that the set $\psi[A]$ has outer measure 1 and inner measure 0 (in
$\prod\limits_{k<\omega}N_k\times\can$), so we may apply $(\boxtimes)_{\rm
sup}^4$ to it, and we get a suitable function $h:\prod\limits_{k<\omega}N_k
\longrightarrow\can$. Let $B=\{x\in\can:(\exists y)((x,y)\in C)\}$, and let
$f^*:B\longrightarrow\can$ be defined by
\[(x,f^*(x))=\psi^{-1}\big((\psi_0(x,y),h(\psi_0(x,y)))\big)\]
for some (equivalently: all) $y$ such that $(x,y)\in C$. Easily $f^*$ is a
Borel function. Take any Borel extension $f:\can\longrightarrow\can$ of
$f^*$ - it is as required in $(\boxtimes)_{\rm sup}^3$ for $A$. 
\medskip

\noindent $(\boxtimes)_{\rm sup}^3\Rightarrow(\boxtimes)_{\rm sup}^4$:\qquad
Even easier. (Note that, since all $N_k$'s are powers of $2$, we have a very
nice measure preserving homeomorphism $\psi^*:\prod\limits_{k<\omega}N_k
\longrightarrow\can$.)
\end{proof}

\begin{theorem}
\label{main}
It is consistent that every sup-measurable function is Lebesgue measurable.
\end{theorem}

\begin{proof}
Start with universe $\bV$ satisfying CH. Let $\bar{\bQ}=\langle\bP_\alpha,
\dot{\bQ}_\alpha:\alpha<\omega_2\rangle$ be countable support iteration such
that each iterand $\dot{\bQ}_\alpha$ is (forced to be) the forcing notion
$\starQ$ (defined in the second section; of course it is taken in the
respective universe $\bV^{\bP_\alpha}$). It follows from \ref{corprobou}
(and \cite[Ch.~VI, 2.8D]{Sh:f}) that the limit $\bP_{\omega_2}$ is proper
and $\baire$--bounding. Also it satisfies $\aleph_2$--cc, and consequently
the forcing with $\bP_{\omega_2}$ does not collapse cardinals nor changes
cofinalities (and $\forces_{\bP_{\omega_2}}\mbox{`` }\con=\aleph_2\mbox{
''}$).  

We are going to prove that
\[\forces_{\bP_{\omega_2}}\mbox{`` every sup-measurable function is Lebesgue 
measurable ''}.\]
By \ref{boxtimes}, it is enough to show that $\forces_{\bP_{\omega_2}}
(\boxtimes)^4_{\rm sup}$. To this end suppose that $\dot{A}$ is a
$\bP_{\omega_2}$--name for a subset of $\prod\limits_{k<\omega}N_k\times
\can$ such that both $\dot{A}$ and its complement are of outer measure
one. By a standard argument using $\aleph_2$--cc of $\bP_{\omega_2}$ (and
the fact that each $\bP_\alpha$ for $\alpha<\omega_2$ has a dense subset of
size $\aleph_1$), we may find $\delta<\omega_2$ of cofinality $\omega_1$,
and a $\bP_\delta$--name $\dot{A}_\delta$ such that  
\[\begin{array}{ll}
\forces_{\bP_{\omega_2}}&\mbox{`` }\dot{A}\cap(\prod\limits_{k<\omega}N_k
\times\can)^{\bV^{\bP_\delta}}=\dot{A}_\delta\mbox{ ''},\qquad\mbox{ and}\\
\forces_{\bP_\delta}&\mbox{`` }\dot{A}_\delta\mbox{ has outer measure $1$
and inner measure $0$ ''}.
  \end{array}\]
Let $\dot{h}$ be the $\bP_{\delta+1}$--name for the continuous function from
$\prod\limits_{k<\omega}N_k$ to $\can$ added at stage $\delta+1$ by
$\dot{\bQ}_\delta=(\starQ)^{\bV^{\bP_\delta}}$ (as defined right before
\ref{crucial}). Then, by \ref{crucial} (applied to $\dot{A}_\delta$ and to
its complement), in $\bV^{\bP_{\delta+1}}$ the set 
\[X_\delta\stackrel{\rm def}{=}\{x\in\prod\limits_{k<\omega} N_k:(x,
\dot{h}(x))\in\dot{A}_\delta\}\]
has outer measure 1 and inner measure 0. By \ref{snep} + \ref{nullpres}, in
$\bV^{\bP_{\delta+1}}$ we may use \cite[Thm 7.8]{Sh:630} (and
\cite[Ch.~XVIII, 3.8]{Sh:f}) to conclude that $\bP_{\omega_2}/
\bP_{\delta+1}$ preserves non-nullity of sets from $\bV^{\bP_{\delta+1}}$. 
Consequently, 
\[\forces_{\bP_{\omega_2}}\mbox{`` the set $X_\delta$ and its complement
have outer measure one ''},\]
finishing the proof.  
\end{proof}

\begin{remark}
\label{remark}
Note that for the iteration $\langle\bP_\alpha,\dot{\bQ}_\alpha:\alpha<
\omega_2\rangle$ to work for the proof of \ref{main} we do not need that all
iterands are $\starQ$. It is enough that for some stationary set $S\subseteq
\{\delta<\omega_2:\cf(\delta)=\omega_1\}$, for every $\alpha\in S$, we have
$\forces_{\bP_\alpha}\dot{\bQ}_\alpha=\starQ$, and that the forcings used in
the iteration are such that each $\bP_{\omega_2}/\bP_{\delta+1}$ preserves
non-nullity of sets from $\bV^{\bP_{\delta+1}}$. So, in particular, we may
use in the iteration also other (s)nep forcing notions preserving ``old
reals are not null''. This will be used in the next section, where we will 
add the random forcing ``here-and-there''.
\end{remark}

\section{Possibly every real function is continuous on a non-null set}
The aim of this section is to show that a slight modification of the
iteration from the previous section results in a model in which every
function $f:\mbR\longrightarrow\mbR$ agrees with a continuous function on a
set of positive outer measure. Let us start with a reduction of the problem
that exposes relevance of the tools developed earlier.

\begin{proposition}
\label{reduction}
Assume:
\begin{enumerate}
\item[(a)] the condition $(\boxtimes)_{\rm sup}^3$ of \ref{boxtimes} holds
true,
\item[(b)] for every function $f^*:\can\longrightarrow\can$ there are
functions $f_1,f_2$ and a set $A$ such that  
\begin{itemize}
\item $A\subseteq\can$ and $f_1:A\longrightarrow\can$ is such that the set 
\[\{(x,f_1(x)):x\in A\}\subseteq\can\times\can\]
has positive outer measure,   
\item $f_2:\can\times\can\longrightarrow\can$ is Borel, and
\item $(\forall x\in A)(f^*(x)=f_2(x,f_1(x)))$.
\end{itemize}
\end{enumerate}
Then for every function $f:\mbR\longrightarrow\mbR$ there is a continuous 
function $g:\mbR\longrightarrow\mbR$ such that the set $\{x\in\mbR:f(x)
=g(x)\}$ has positive outer measure. 
\end{proposition}

\begin{proof}
Assume $f:\mbR\longrightarrow\mbR$. Let $\varphi:\mbR\longrightarrow\can$ be
a Borel isomorphism preserving null sets (see, e.g., \cite[Thm
17.41]{Ke94}), and let $f^*=\varphi\comp f\comp\varphi^{-1}$. Let
$f_1,f_2,A$ be given by the assumption (b) for $f^*$. Put
$A^*=\{(x,f_1(x)):x\in A\}\subseteq\can\times\can$. We know that $A^*$ is a
non-null set (and consequently it is non-measurable), so applying
$(\boxtimes)^3_{\rm sup}$ we may pick a Borel function $g_0:\can
\longrightarrow\can$ such that the set  
\[B\stackrel{\rm def}{=}\{x\in A: f_1(x)= g_0(x)\}\]
has positive outer measure, and so does $\varphi^{-1}[B]$. Let $g_1:\mbR
\longrightarrow\mbR$ be defined by 
\[g_1(x)=\varphi^{-1}\Big(f_2\big(\varphi(x),g_0(\varphi(x))\big)\Big).\]
Clearly $g_1$ is Borel and for each $x\in\varphi^{-1}[B]$ we have
$g_1(x)=f(x)$. Finally, using Lusin's theorem (see, e.g., \cite[Thm
17.12]{Ke94}) we may pick a continuous function $g:\mbR\longrightarrow\mbR$
such that the set $\{x\in \varphi^{-1}[B]: g_1(x)=g(x)\}$ is not null (just
take $g$ so that it agrees with $g_1$ on a set of large enough measure). 
\end{proof}

The iteration of \ref{main} will be changed by adding random reals on a
stationary set. So just for uniformity of our notation we represent the
random real forcing as $\rQ$. Let $\bH^r(i)=2$ (for $i<\omega$). Let $K^r$
consist of tree creatures $t\in\TCR[\bH^r]$ 
such that
\begin{itemize}
\item $\dis[t]=(k_t,\eta_t,P_t)$, where $k_t<\omega$, $\eta_t\in
\prod\limits_{i<k_t}\bH^r(i)$, $\emptyset\neq P_t\subseteq 2$, and 
\item $\nor[t]=k_t$,
\item $\val[t]=\{\langle\eta_t,\nu\rangle:\eta_t\vartriangleleft\nu\in
\prod\limits_{i\leq k_t}\bH^r(i)\ \&\ \nu(k_t)\in P_t\}$.
\end{itemize}
The operation $\Sigma^r$ is trivial:
\[\Sigma^r(t)=\{s\in K^r:\eta_s=\eta_t\ \&\ P_s\subseteq P_t\}.\]
For $t\in K^r$ and a sequence $\langle r_\nu:\nu\in\pos(t)\rangle\subseteq
[0,1]$ we let 
\[F^r_t(r_\nu:\nu\in\pos(t))=\frac{\sum\{r_\nu:\nu\in\pos(t)\}}{2}.\]
It is easy to check that $(K^r,\Sigma^r,\bF^r)$ is a (nice) measured tree
creating triple for $\bH^r$, and that the forcing notion $\rQ$ is
(equivalent to) the random real forcing.
\medskip

Like in \ref{main}, we start with universe $\bV$ satisfying CH. Let
$Z\subseteq \{\delta<\omega_2:\cf(\delta)=\omega_1\}$ be a stationary set
such that $\{\delta<\omega_2:\cf(\delta)=\omega_1\}\setminus Z$ is
stationary as well. Let $\bar{\bQ}=\langle\bP_\alpha,\dot{\bQ}_\alpha:\alpha
<\omega_2\rangle$ be countable support iteration such that
\begin{itemize}
\item if $\alpha\in Z$, then $\forces_{\bP_\alpha}\dot{\bQ}_\alpha=\rQ$,
\item if $\alpha\in\omega_2\setminus Z$, then $\forces_{\bP_\alpha}
\dot{\bQ}_\alpha=\starQ$.
\end{itemize}
We are going to show that 
\[\forces_{\bP_{\omega_2}}\mbox{`` every real function is continuous on a
non-null set ''},\]
and for this we will show that the assumptions of \ref{reduction} are
satisfied in $\bV^{\bP_{\omega_2}}$. First note that $\bP_{\omega_2}\forces
(\boxtimes)^3_{\rm sup}$ (see \ref{remark}; remember \ref{boxtimes}). To
show that, in $\bV^{\bP_{\omega_2}}$, the assumption (b) of \ref{reduction}
holds, we need to analyze conditions and continuous reading of names in the
iteration.  

\begin{definition}
\label{gfc}
Let $(K,\Sigma,\bF)$ be a measured tree creating triple for $\bH$ (say,
either $(K^*,\Sigma^*,\bF^*)$ defined in the second section, or $(K^r,
\Sigma^r,\bF^r)$ defined above). 
\begin{enumerate}
\item {\em A finite candidate for $(K,\Sigma,\bF)$} (or just for $(K,
\Sigma)$) is a system $\bs=\langle s_\eta:\eta\in S\setminus\max(S)\rangle$
such that   
\begin{itemize}
\item $S\subseteq\bigcup\limits_{n<\omega}\prod\limits_{i<n}\bH(i)$ is a
finite tree, $s_\eta\in K\cap\TCR_\eta[\bH]$ for $\eta\in S\setminus\max(
S)$, 
\item $\max(S)\subseteq\prod\limits_{i<m}\bH(i)$ for some $m=\rht(\bs)$ (we
will call this $m$ the height of the candidate $\bs$),
\item if $\eta\in S\setminus\max(S)$, then $\suc_S(\eta)=\pos(s_\eta)$. 
\end{itemize}
We may also write $\mrot(\bs)$ for $\mrot(S)$ (and call it the root of the
candidate $\bs$), and write $\max(\bs)$ for $\max(S)$. 
\item Let $\FC(K,\Sigma)$ be the family of all finite candidates for $(K,
\Sigma)$.  
\item For candidates $\bs^0,\bs^1\in\FC(K,\Sigma)$, we say that {\em $\bs^1$
end--extends $\bs^0$} (in short: $\bs^0\preceq_{\rm end}\bs^1$) if
$\mrot(\bs^1)=\mrot(\bs^0)$, $\rht(\bs^1)\geq\rht(\bs^0)$ and, letting 
$\bs^\ell=\langle s^\ell_\eta:\eta\in S^\ell\setminus\max(S^\ell)\rangle$,
we have $S^0\subseteq S^1$ and $(\forall\eta\in S^0\setminus\max(S^0))(
s^0_\eta=s^1_\eta)$. 
\item We say that a condition $p\in\bQ^\tree_\emptyset(K,\Sigma)$ {\em
end--extends\/} a candidate $\bs=\langle s_\eta:\eta\in S\setminus\max(S)
\rangle\in\FC(K,\Sigma)$ if 
\[\mrot(p)=\mrot(\bs),\quad S\subseteq T^p\quad\mbox{ and }\quad s_\eta=
t^p_\eta\ \mbox{ for }\ \eta\in S\setminus\max(S).\]  
\end{enumerate}
\end{definition}

\begin{definition}
\label{pretemplate}
\begin{enumerate}
\item {\em A finite pre--template\/} is a tuple 
\[\bt=\langle w^\bt,\bk^\bt,\bc^\bt,\bar{\cY}^\bt\rangle=\langle w,\bk,\bc,
\bar{\cY}\rangle\] 
such that 
\begin{enumerate}
\item[$(\alpha)$] $w$ is a finite non-empty set of ordinals below
$\omega_2$, $w=\{\alpha_0,\ldots,\alpha_n\}$ (the increasing enumeration);

let $x_i$ be $r$ if $\alpha_i\in Z$, and $x_i$ be $*$ if $\alpha_i\in
\omega_2\setminus Z$,
\item[$(\beta)$]  $\bk:w\longrightarrow\omega$, $\bc=\langle c_{\alpha_0},
\ldots,c_{\alpha_n}\rangle$, $\bar{\cY}=\langle \cY_{\alpha_0},\ldots,
\cY_{\alpha_n} 
\rangle$ (we treat $\bc,\bar{\cY}$ as functions with domain $w$),
\item[$(\gamma)$] $c_{\alpha_0}\in\FC(K^{x_0},\Sigma^{x_0})$, $\rht(
c_{\alpha_0})=\bk(\alpha_0)$, $\cY_{\alpha_0}=\{\langle s\rangle: s\in
\max(c_{\alpha_0})\}$, and for $0<i\leq n$: 
\item[$(\delta)$] $c_{\alpha_i}:\cY_{\alpha_{i-1}}\longrightarrow\FC(
K^{x_i},\Sigma^{x_i})$ is such that $\rht(c_{\alpha_i}(\bar{\nu}))=\bk(
\alpha_i)$ for each $\bar{\nu}\in\cY_{\alpha_{i-1}}$,\\
$\cY_{\alpha_i}=\{\bar{\nu}\conc\langle\nu_{\alpha_i}\rangle:\bar{\nu}=
\langle\nu_{\alpha_0},\ldots,\nu_{\alpha_{i-1}}\rangle\in\cY_{\alpha_{i-1}}\ 
\&\ \nu_{\alpha_i}\in\max(c_{\alpha_i}(\bar{\nu}))\}$.  

(We think of elements of $\cY_{\alpha_i}$ as functions from $\{\alpha_0,
\ldots, \alpha_i\}$ with values being sequences in appropriate
$\prod\limits_{j<\bk(\alpha_\ell)}\bH^{x_\ell}(j)$.)
\end{enumerate}
$\cY_{\alpha_n}$ will be also called $\cY_*$ or $\cY^\bt_*$.
\item We say that a finite pre--template $\bt'$ {\em properly extends\/} a
pre--template $\bt$ (and then we write $\bt\preceq\bt'$) if
\begin{enumerate}
\item[$(\alpha)$] $w^\bt\subseteq w^{\bt'}$, and $(\forall\alpha\in w^\bt)(
\bk^\bt(\alpha)\leq \bk^{\bt'}(\alpha))$, and
\item[$(\beta)$] let $w^{\bt'}=\{\alpha_0,\ldots,\alpha_n\}$ (the
increasing enumeration).

If $\ell^*=\min\{i\leq n:\alpha_i\in w^\bt\}$, then for every $\langle
\nu_{\alpha_0},\ldots,\nu_{\alpha_{\ell^*-1}}\rangle\in\cY^{\bt'}_{
\alpha_{\ell^*-1}}$ we have $c^\bt_{\alpha_{\ell^*}}\preceq_{\rm end}
c^{\bt'}_{\alpha_{\ell^*}}(\nu_{\alpha_0},\ldots,\nu_{\alpha_{\ell^*-1}})$. 

If $\ell>\ell^*$ is such that $\alpha_\ell\in w^\bt$ and $k<\ell$ is such
that $\alpha_k$ is the predecessor of $\alpha_\ell$ in $w^\bt$, then for
every $\langle\nu_{\alpha_0},\ldots,\nu_{\alpha_{\ell-1}}\rangle\in
\cY^{\bt'}_{\alpha_{\ell-1}}$ we have 
\[\begin{array}{l}
\langle\nu_{\alpha_i}\restriction\bk^\bt(\alpha_i):i<\ell\ \&\ \alpha_i\in
w^\bt \rangle\in\cY^\bt_{\alpha_k}\quad\mbox{ and}\\
c^\bt_{\alpha_\ell}(\nu_{\alpha_i}\restriction\bk^\bt(\alpha_i):i<\ell\ \&\
\alpha_i\in w^\bt)\preceq_{\rm end} c^{\bt'}_{\alpha_\ell}(\nu_{\alpha_0},
\ldots,\nu_{\alpha_{\ell-1}}).
\end{array}\]
\end{enumerate}
\item For an ordinal $\zeta<\omega_2$ and a finite pre-template $\bt$ we
define the restriction $\bt'=\bt\restriction\zeta$ of $\bt$ in a natural
way: $w^{\bt'}=w^\bt\cap\zeta$, $\bk^{\bt'}=\bk^\bt\restriction w^{\bt'}$,
$\bc^{\bt'}=\bc^\bt\restriction w^{\bt'}$ and $\bar{\cY}^{\bt'}=
\bar{\cY}^\bt \restriction w^{\bt'}$. (Note that $\bt\restriction\zeta
\preceq\bt$.)
\item We say that finite pre-templates $\bt,\bt'$ are {\em isomorphic\/} if
$|w^\bt|=|w^{\bt'}|$, and if $h:w^\bt\longrightarrow w^{\bt'}$ is the order
preserving isomorphism, then
\begin{itemize}
\item $h[w^\bt\cap Z]=w^{\bt'}\cap Z$, and 
\item $\bk^\bt=\bk^{\bt'}\comp h$,  $\bc^\bt=\bc^{\bt'}\comp h$, and 
$\bar{\cY}^\bt=\bar{\cY}^{\bt'}\comp h$.
\end{itemize}
We also may say that {\em $h$ is an isomorphism from $\bt$ to $\bt'$}.
\end{enumerate}
\end{definition}

\begin{definition}
\label{obeys}
By induction on $n=|w^\bt|-1$ we define 
\begin{enumerate}
\item[(a)] when a condition $p\in\bP_{\omega_2}$ obeys a pre-template $\bt$,
and 
\item[(b)] if $w^\bt=\{\alpha_0,\ldots,\alpha_n\}$, $\bar{\nu}=\langle
\nu_{\alpha_0},\ldots,\nu_{\alpha_n}\rangle\in\cY^\bt_*$, and $p\in
\bP_{\omega_2}$ obeys $\bt$, then we define a condition
$p^{[\bt,\bar{\nu}]}\in\bP_{\omega_2}$ stronger that $p$.
\end{enumerate}
First consider the case when $n=0$. Let $\bt$ be a pre-template such that
$w^\bt=\{\alpha_0\}$ and let $p\in\bP_{\omega_2}$. We say that $p$ obeys
$\bt$ if 
\[p\restriction\alpha_0\forces_{\bP_{\alpha_0}}\mbox{`` }p(\alpha_0)\mbox{
end extends the candidate }c^\bt_{\alpha_0}\mbox{ ''.}\]
If $p$ obeys $\bt$ as above, and $\bar{\nu}=\langle\nu_{\alpha_0}\rangle
\in\cY^\bt_{\alpha_0}$, then $p^{[\bt,\bar{\nu}]}$ is defined as follows: 
\begin{itemize}
\item $p^{[\bt,\bar{\nu}]}\restriction(\omega_2\setminus\{\alpha_0\})=
p\restriction(\omega_2\setminus\{\alpha_0\})$,\quad and
\item $p^{[\bt,\bar{\nu}]}\restriction\alpha_0\forces_{\bP_{\alpha_0}}$``
$p^{[\bt,\bar{\nu}]}(\alpha_0)=(p(\alpha_0))^{[\nu_{\alpha_0}]}$ ''.
\end{itemize}
Now, suppose that $w^\bt=\{\alpha_0,\ldots,\alpha_n\}$ (the increasing
enumeration; $n>0$), and that we have dealt with $n-1$ already. We say that
a condition $p\in\bP_{\omega_2}$ obeys $\bt$ if  
\begin{itemize}
\item $p$ obeys $\bt\restriction\alpha_n$, and
\item for every $\bar{\nu}=\langle\nu_{\alpha_0},\ldots,\nu_{\alpha_{n-1}}
\rangle\in \cY^\bt_{\alpha_{n-1}}$, the condition $p^{[\bt\restriction
\alpha_n,\bar{\nu}]}\restriction\alpha_n$ forces (in $\bP_{\alpha_n}$) that 
$p(\alpha_n)$ end--extends the candidate $c^\bt_{\alpha_n}(\bar{\nu})$.
\end{itemize}
In that case we also define $p^{[\bt,\bar{\nu}]}$ for $\bar{\nu}=\langle
\nu_{\alpha_0},\ldots,\nu_{\alpha_n}\rangle\in\cY^\bt_{\alpha_n}$: 
\begin{itemize}
\item $p^{[\bt,\bar{\nu}]}\restriction\omega_2\setminus\{\alpha_n\}= 
p^{[\bt\restriction\alpha_n,\bar{\nu}\restriction\alpha_n]}\restriction
\omega_2\setminus\{\alpha_n\}$, 
\item $p^{[\bt,\bar{\nu}]}\restriction\alpha_n\forces_{\bP_{\alpha_n}}$``
$p^{[\bt,\bar{\nu}]}(\alpha_n)=(p(\alpha_n))^{[\nu_{\alpha_n}]}$ ''.
\end{itemize}
\end{definition}

\begin{definition}
\label{weaktemplate}
\begin{enumerate}
\item {\em A weak template\/} is a $\preceq$--increasing sequence
$\bar{\bt}=\langle\bt_n:n<\omega\rangle$ of finite pre-templates such that  
\[(\forall\alpha\in\bigcup_{n<\omega}w^{\bt_n})(\lim_{n\to\infty}\bk^{\bt_n}
(\alpha)=\infty).\]
\item We say that {\em weak templates $\bar{\bt},\bar{\bt}'$ are
isomorphic\/} if  
\begin{itemize}
\item $\otp(\bigcup\limits_{n<\omega} w^{\bt_n})=\otp(\bigcup\limits_{n<
\omega} w^{\bt_n'})$, and
\item letting $h:\bigcup\limits_{n<\omega} w^{\bt_n}\longrightarrow
\bigcup\limits_{n<\omega} w^{\bt_n'}$ be the order isomorphism, we have that
all restrictions $h\restriction w^{\bt_n}$ (for $n<\omega$) are isomorphisms
from $\bt_n$ to $\bt_n'$.  
\end{itemize}
(We will also call the mapping $h$ as above {\em the isomorphism from
$\bar{\bt}$ to $\bar{\bt}'$\/}.) 
\item A condition $p\in\bP_{\omega_2}$ {\em obeys the weak template
$\bar{\bt}=\langle\bt_n:n<\omega\rangle$} if $\dom(p)=\bigcup\limits_{n<
\omega}w^{\bt_n}$ and $p$ obeys each pre-template $\bt_n$ (for $n<\omega$).
\item {\em A weak template with a name\/} is a pair $(\bar{\bt},\bar{\tau})$
such that $\bar{\bt}=\langle\bt_n:n<\omega\rangle$ is a weak template, and
$\bar{\tau}=\langle\tau_n:n<\omega\rangle$ is a sequence of functions such
that $\tau_n:\cY^{\bt_n}_*\longrightarrow 2^{\textstyle n}$, and if 
$\langle\nu_\alpha:\alpha\in w^{\bt_{n+1}}\rangle\in\cY^{\bt_{n+1}}_*$, then  
\[\tau_n(\nu_\alpha\restriction\bk^{\bt_n}(\alpha):\alpha\in w^{\bt_n})
\vartriangleleft \tau_{n+1}(\nu_\alpha:\alpha\in w^{\bt_{n+1}}).\]
\item Let $(\bar{\bt},\bar{\tau}),(\bar{\bt}',\bar{\tau}')$ be weak
templates with names. We say that they are isomorphic provided that $\bt$
and $\bt'$ are isomorphic, and the isomorphism maps $\bar{\tau}$ to
$\bar{\tau}'$. (To be more precise, if $h$ is the isomorphism from $\bt$ to 
$\bt'$, then for each $n<\omega$ it induces a bijection $g_n:\cY^{\bt_n}_*
\longrightarrow\cY^{\bt_n'}_*$; we request that $\tau_n=\tau_n'\comp g_n$.)  
\item Let $(\bar{\bt},\bar{\tau})$ be a weak template with a name, $p\in
\bP_{\omega_2}$ and let $\dot{\tau}$ be a $\bP_{\omega_2}$--name for a real
in $\can$. We say that {\em $(p,\dot{\tau})$ obeys $(\bar{\bt},\bar{\tau})$}
if 
\begin{itemize}
\item the condition $p$ obeys the weak template $\bt$, and 
\item for each $n<\omega$ and $\bar{\nu}\in\cY^{\bt_n}_*$ we have:
$p^{[\bt_n,\bar{\nu}]}\forces_{\bP_{\omega_2}}\dot{\tau}\restriction n=
\tau_n(\bar{\nu})$.
\end{itemize}
\end{enumerate}
\end{definition}

\begin{lemma}
\label{easycount}
\begin{enumerate}
\item There are only countably many isomorphism types of finite pre-templates.
\item There are $\con$ many isomorphism types of weak templates with names. 
\end{enumerate}
\end{lemma}

\begin{lemma}
\label{condec}
Suppose that $\dot{\tau}$ is a $\bP_{\omega_2}$--name for a real in $\can$
and $p\in\bP_{\omega_2}$. Then there is a condition $q\in\bP_{\omega_2}$
stronger than $p$, and a weak template with a name $(\bt,\bar{\tau})$ such
that $(q,\dot{\tau})$ obeys $(\bt,\bar{\tau})$.
\end{lemma}

\begin{proof}
It is a standard application of fusion arguments somewhat similar to the
proof of Baumgartner \cite[Lemma 7.3]{B3}; compare also the proof of
preservation of ``proper+$\baire$--bounding'' in \cite[Ch. VI]{Sh:f} or
\cite{Go}. (Of course, we use here Lemma \ref{preconrea}. 
\end{proof}

Note that there are weak templates $\bt$ such that no condition $p\in
\bP_{\omega_2}$ obeys $\bt$ -- there could be a problem with norms and/or
measures! From all weak templates we will select only those which correspond
to conditions in $\bP_{\omega_2}$ (and they will be called just {\em
templates}; see \ref{templates} below). 

\begin{definition}
\label{covers}
\begin{enumerate}
\item {\em A cover for a condition\/} $p\in\bQ^\tree_\emptyset(K^*,
\Sigma^*)$ is a condition $q\in\bQ^\tree_\emptyset(K^*,\Sigma^*)$ defined so
that $\mrot(p)=\mrot(q)$, $q\leq p$ and:\\ 
if $\eta\in T^p$, $k=\lh(\eta)$, then $\nor[t^q_\eta]=\nor[t^p_\eta]$,
$g_{t^q_\eta}=g_{t^p_\eta}$, and 
\[P_{t^q_\eta}=\{f\in\bH^*(k): g_{t^q_\eta}\subseteq f\},\]
if $\eta\notin T^p$, $k=\lh(\eta)$, then $g_{t^q_\eta}=\emptyset$,
$P_{t^q_\eta}=\bH^*(k)$ and $\nor[t^q_\eta]=k$.
\item Let $p\in \bQ^\tree_\emptyset(K^*,\Sigma^*)$, and let $q$ be the cover
of $p$, and assume that $T^q$ is a perfect tree. {\em The covering mapping
for $p$\/} is the natural homeomorphism $\bh_p:[T^q]\longrightarrow\can$,
defined as follows. First we define a mapping $h_p:T^q\longrightarrow\fs$:
we let $h_p(\mrot(T^q))=\langle\rangle$. Suppose that $h_p(\eta)$ has been
defined, $\eta\in T^q$, and say $h_p(\eta)\in 2^{\textstyle n}$,
$n<\omega$. We note that $|\pos(t^q_\eta)|$ is a power of $2$, and thus we
may pick $k>n$ such that $|\pos(t^q_\eta)|=|2^{\textstyle [n,k)}|$. Now,
$h_p$ maps $\pos(t^q_\eta)$ onto $\{\nu\in 2^{\textstyle
k}:h_p(t^q_\eta)\vartriangleleft\nu\}$ (preserving some fixed well-ordering
of ${\mathcal H}(\aleph_1)$). Finally we let $\bh_p(\rho)=
\bigcup\limits_{n<\omega}h_p(\rho\restriction n)$.   

In the case when $T^q$ is not a perfect tree, then the covering mapping for
$p$ maps all elements of $[T^q]$ to the sequence with constant value $0$. 

\item Similarly we define {\em the covering mapping\/} $\bh_p$ for a
condition $p\in\bQ^\tree_\emptyset(K^r,\Sigma^r)$ such that $T^p$ is a
perfect tree. So, first we let $h_p:T^p\longrightarrow\fs$ be a
$\trianglelefteq$--preserving mapping such that $h_p(\mrot(T^p))=\langle
\rangle$, if $\eta\in T^p$ is a splitting node of $T^p$ then
$h_p(\eta\conc\langle\ell\rangle)=h_p(\eta)\conc\langle\ell\rangle$ (for
$\ell<2$), $\rng(h_p)=\fs$. Next, $\bh_p$ is the homeomorphism from $[T^p]$
onto $\can$ induced by $h_p$.  (And if $T^p$ is not a perfect tree, then the
covering mapping $\bh_p$ is constant.)

In this case also {\em the cover for the condition $p$} is $p$ itself.
\end{enumerate}
\end{definition}

Now we are going to introduce the main technical tool involved in the proof
that our iteration is OK. Fix a weak template $\bar{\bt}=\langle\bt_n:n<
\omega\rangle$ for a while. Let $w^{\bar{\bt}}=\bigcup\limits_{n<\omega}
w^{\bt_n}$ and $\zeta_{\bar{\bt}}=\otp(w^{\bar{\bt}})$, and let
$w^{\bar{\bt}}=\langle\alpha_\zeta:\zeta<\zeta_{\bar{\bt}}\rangle$ (the
increasing enumeration). For $\zeta<\zeta_{\bar{\bt}}$ let $x_\zeta$ be $r$
if $\alpha_\zeta\in Z$, and $*$ if $\alpha_\zeta\notin Z$. 

By induction on $\zeta\leq\zeta_{\bar{\bt}}$ we define a space
$\cZ^{\bar{\bt}}_\zeta$ and mappings
\[\pi^{\bar{\bt}}_\zeta:\cZ^{\bar{\bt}}_\zeta\longrightarrow
\bQ^\tree_\emptyset(K^{x_\zeta},\Sigma^{x_\zeta})\qquad\mbox{ and }\qquad
\psi^{\bar{\bt}}_\zeta:\cZ^{\bar{\bt}}_\zeta\longrightarrow
(\can)^{\textstyle\zeta}.\] 
First we let $\cZ^{\bar{\bt}}_0=\{\emptyset\}$ and $\pi^{\bar{\bt}}_0(
\emptyset)\in\bQ^\tree_\emptyset(K^{x_0},\Sigma^{x_0})$ is a condition
end--extending all $c^{\bt_n}_{\alpha_0}$ (for $n<\omega$, $\alpha_0\in
w^{\bt_n}$) (and $\psi^{\bar{\bt}}_\zeta(\emptyset)=\emptyset$). 

\noindent Suppose now that $\zeta+1\leq \zeta_{\bar{\bt}}$ and we have
defined $\cZ^{\bar{\bt}}_\zeta,\pi^{\bar{\bt}}_\zeta$ and
$\psi^{\bar{\bt}}_\zeta$. We let 
\[\cZ^{\bar{\bt}}_{\zeta+1}=\{\bar{z}\conc\langle z_\zeta\rangle:\bar{z}\in
\cZ^{\bar{\bt}}_\zeta\ \&\ z_\zeta\in [T^{\pi^{\bar{\bt}}_\zeta(\bar{z})}]
\subseteq \prod_{i<\omega}\bH^{x_\zeta}(i)\},\]
and let $\bar{z}^*=\langle z_0,\ldots,z_\zeta\rangle=\bar{z}\conc\langle
z_\zeta\rangle\in\cZ^{\bar{\bt}}_{\zeta+1}$ (we ignore the first term
``$\emptyset$'' of the sequence $\bar{z}$). To define $\psi^{\bar{\bt}}_{ 
\zeta+1}(\bar{z}^*)$, we let $q\in\bQ^\tree_\emptyset(K^{x_\zeta},\Sigma^{
x_\zeta})$ be a cover for the condition $\pi^{\bar{\bt}}_\zeta(\bar{z})$,
and let $\bh:[T^q]\longrightarrow\can$ be the covering mapping for
$\pi^{\bar{\bt}}_\zeta(\bar{z})$ (see \ref{covers}). Put $\psi^{\bar{
\bt}}_{\zeta+1}(\bar{z}^*)=\psi^{\bar{\bt}}_\zeta(\bar{z})\conc\langle\bh(
z_\zeta)\rangle$.

\noindent If $\zeta+1<\zeta_{\bar{\bt}}$, then we also define $\pi^{\bar{
\bt}}_{\zeta+1}(\bar{z}^*)$ as a condition in $\bQ^\tree_\emptyset(
K^{x_{\zeta+1}},\Sigma^{x_{\zeta+1}})$ such that 
\begin{itemize}
\item if $n<\omega$, $w^{\bt_n}=\{\alpha_{\zeta_0},\ldots,\alpha_{\zeta_m}
\}$ (the increasing enumeration), and $\zeta_\ell=\zeta+1$, $\ell\leq m$,
{\em then\/} $\pi^{\bar{\bt}}_{\zeta+1}(\bar{z}^*)$ end extends 
$c^{\bt_n}_{\alpha_{\zeta_\ell}}(z_{\zeta_0}\restriction \bk^{\bt_n}(
\alpha_{\zeta_0}),\ldots, z_{\zeta_{\ell-1}}\restriction \bk^{\bt_n}(
\alpha_{\zeta_{\ell-1}}))$. 
\end{itemize}

\noindent Suppose now that $\zeta\leq\zeta_{\bar{\bt}}$ is a limit ordinal,
and that we have defined $\cZ^{\bar{\bt}}_\xi$, $\pi^{\bar{\bt}}_\xi$ and
$\psi^{\bar{\bt}}_\xi$ for $\xi<\zeta$. We put 
\[\cZ^{\bar{\bt}}_\zeta=\{\langle z_\rho:\rho<\zeta\rangle: (\forall\xi<
\zeta)(\langle z_\rho:\rho<\xi\rangle\in \cZ^{\bar{\bt}}_\xi)\}\]
(again, above, like before and later, we ignore the first term
``$\emptyset$'' whenever considering elements of $\cZ^{\bar{\bt}}_\xi$). 
The mapping $\psi^{\bar{\bt}}_\zeta:\cZ^{\bar{\bt}}_\zeta\longrightarrow
(\can)^{\textstyle\zeta}$ is such that $\psi^{\bar{\bt}}_\zeta(\bar{z})
\restriction\xi=\psi^{\bar{\bt}}_\xi(\bar{z}\restriction\xi)$ (for
$\bar{z}\in\cZ^{\bar{\bt}}_\zeta$). Also if, additionally, $\zeta<
\zeta_{\bar{\bt}}$, then for $\bar{z}=\langle z_\rho:\rho<\zeta\rangle\in
\cZ^{\bar{\bt}}_\zeta$ we let $\pi^{\bar{\bt}}_\zeta(\bar{z})$ be an element
of $\bQ^\tree_\emptyset(K^{x_\zeta},\Sigma^{x_\zeta})$ such that 
\begin{itemize}
\item if $n<\omega$, $w^{\bt_n}=\{\alpha_{\zeta_0},\ldots,\alpha_{\zeta_m}
\}$ (the increasing enumeration), and $\zeta_\ell=\zeta$, $\ell\leq m$,
{\em then\/} $\pi^{\bar{\bt}}_\zeta(\bar{z})$ end extends $c^{\bt_n}_{
\alpha_{\zeta_\ell}}(z_{\zeta_0}\restriction \bk^{\bt_n}(\alpha_{\zeta_0}),
\ldots,z_{\zeta_{\ell-1}}\restriction\bk^{\bt_n}(\alpha_{\zeta_{\ell-1}}))$.
\end{itemize}

\begin{definition}
\label{templates}
Let $\bar{\bt}$ be a weak template, and $w^{\bar{\bt}}=\langle
\alpha_\zeta:\zeta<\zeta_{\bar{\bt}}\rangle$ be the increasing enumeration.
Also for $\zeta<\zeta_{\bar{\bt}}$ let $x_\zeta$ be $r$ if $\alpha_\zeta\in
Z$, and $*$ if $\alpha_\zeta\notin Z$. We say that $\bar{\bt}$ is {\em a
template\/} if for every $\zeta<\zeta_{\bar{\bt}}$ and $\bar{z}\in
\cZ^{\bar{\bt}}_\zeta$ we have 
\[\pi^{\bar{\bt}}_\zeta(\bar{z})\in\bQ^\mtree_4(K^{x_\zeta},
\Sigma^{x_\zeta},\bF^{x_\zeta}).\]
\end{definition}

\begin{proposition}
\label{easytemp}
\begin{enumerate}
\item Assume that $p\in\bP_{\omega_2}$. Then there are a condition $q\in
\bP_{\omega_2}$ and a template $\bar{\bt}$ such that $q\geq p$ obeys
$\bar{\bt}$, $\omega\leq\zeta_{\bar{\bt}}<\omega_1$, and for some
enumeration $\langle\zeta_n:n<\omega\rangle$ of $\zeta_{\bar{\bt}}$ we have:  
\begin{enumerate}
\item[$(\boxplus)$] for every $n<\omega$ and $\bar{z}\in\cZ^{\bar{\bt}}_{
\zeta_n}$, 
\[\mu^\bF(\pi^{\bar{\bt}}_{\zeta_n}(\bar{z}))\geq \big(1-2^{-n-10}\big),\]
where $\bF$ is suitably $\bF^r$ or $\bF^*$. 
\end{enumerate}
[If a template $\bar{\bt}$ satisfies $(\boxplus)$ for an enumeration
$\bar{\zeta}=\langle\zeta_n:n<\omega\rangle$ of $\zeta_{\bar{\bt}}$, then we
we will say that {\em $\bar{\bt}$ behaves well for $\bar{\zeta}$}.]
\item For every template $\bar{\bt}$, there is a condition $p\in
\bP_{\omega_2}$ which obeys $\bar{\bt}$. 
\end{enumerate}
\end{proposition}

For a countable ordinal $\zeta$, the space $(\can)^{\textstyle \zeta}$ is
equipped with the product measure $m^{\rm Leb}$ of countably many copies
of $\can$. We will use the same notation $m^{\rm Leb}$ for this measure in
various products (and related spaces), hoping that no real confusion is
caused.  

\begin{lemma}
\label{specclos}
Let $\zeta<\omega_1$. Suppose that $C\subseteq (\can)^{\textstyle\zeta}$ is
a closed set of positive Lebesgue measure. Then there is a closed set $C^*
\subseteq C$ of positive Lebesgue measure such that for each $\xi<\zeta$: 
\begin{enumerate}
\item[$(\otimes)^\xi_{C^*}$] for every $\bar{y}\in (\can)^\xi$, the set 
\[(C^*)_{\bar{y}}\stackrel{\rm def}{=}\{\bar{y}'\in (\can)^{\textstyle
[\xi,\zeta)}:\bar{y}\conc\bar{y}'\in C^*\}\]
is either empty or has positive Lebesgue measure (in $(\can)^{\textstyle
[\xi,\zeta)}$).  
\end{enumerate}
\end{lemma}

\begin{proof}
For  a set $X\subseteq (\can)^{\textstyle\zeta}$, $\xi<\zeta$, and $\bar{y}
\in (\can)^{\textstyle\xi}$ we let
\[(X)_{\bar{y}}\stackrel{\rm def}{=}\{\bar{y}'\in (\can)^{\textstyle [\xi,
\zeta)}:\bar{y}\conc\bar{y}'\in X\}.\]
We may assume that $\zeta\geq\omega$ (otherwise the lemma is easier and
actually included in this case). Fix an enumeration $\zeta=\{\zeta_n:n<
\omega\}$ such that $\zeta_0=0$, and let  
\[e_0=2^{-4}\cdot m^{\rm Leb}(C),\qquad e_{n+1}=2^{-6(n+2)^2}\cdot
(e_n)^{n+2}.\] 
We are going to define inductively a decreasing sequence $\langle C_n:n<
\omega\rangle$ of closed (non-empty) subsets of $C$ such that $C_0=C_1=C$
and  
\begin{enumerate}
\item[$(\circledast)$] for each $m<n$ and $\bar{y}\in (\can)^{\textstyle
\zeta_m}$ we have 
\[\mbox{either}\quad (C_n)_{\bar{y}}=\emptyset\qquad\mbox{ or}\quad m^{\rm
Leb}((C_n)_{\bar{y}})\geq e_m\cdot\Big(1-\sum_{\ell=m+2}^n 4^{-\ell}\Big).\] 
\end{enumerate}
(Note that $(\circledast)$ implies $m^{\rm Leb}(C_n)\geq e_0\cdot
(1-\sum\limits_{\ell=2}^n 4^{-\ell})$.)

Suppose that $C_n$ has been defined already, $n\geq 1$. Let $\{\xi_\ell:\ell
\leq \ell^*\}$ enumerate the set
\[\{\zeta_m:m\leq n\ \&\ \zeta_m\leq\zeta_n\}\]
in the increasing order. By downward induction on $0<\ell\leq\ell^*$ we
chose open sets $U_\ell\subseteq (\can)^{\textstyle\xi_\ell}$. So, the set
$U_{\ell^*}\subseteq (\can)^{\textstyle\xi_{\ell^*}}$ is such that (remember
$\xi_{\ell^*}=\zeta_n$): 
\begin{itemize}
\item $(\forall\bar{y}\in (\can)^{\textstyle\zeta_n}\setminus U_{\ell^*})
\big(m^{\rm Leb}((C_n)_{\bar{y}})\geq e_n\big)$, 
\item $m^{\rm Leb}\big(C_n\cap (U_{\ell^*}\times (\can)^{\textstyle[\zeta_n,
\zeta)})\big)<e_n$. 
\end{itemize}
Now suppose that $U_{\ell^*},\ldots,U_{\ell+1}$ have been chosen already so
that  
\[m^{\rm Leb}\big(C_n\cap(U_k\times (\can)^{\textstyle [\xi_k,\zeta)})\big)
<\Big(\frac{2^{3n+3}}{e_{n-1}}\Big)^{\ell^*-k}\cdot e_n\]
for each $k\in\{\ell+1,\ldots,\ell^*\}$. Let
\[U=U_{\ell+1}\times (\can)^{\textstyle [\xi_{\ell+1},\zeta)}\cup\ldots\cup
U_{\ell^*}\times (\can)^{\textstyle [\xi_{\ell^*},\zeta)}.\]
Note that (by our assumptions)
\[m^{\rm Leb}(U\cap C_n)<(\ell^*-\ell)\cdot\Big(\frac{2^{3n+3}}{e_{n-1}}\Big 
)^{\ell^*-\ell-1}\cdot e_n.\]
Let $\xi_\ell=\zeta_m$ and $A=\{\bar{y}\in(\can)^{\textstyle\zeta_m}: m^{\rm
Leb}\big((C_n\cap U)_{\bar{y}}\big)>\frac{e_m}{2^{2n+2}}\}$. Note that 
\[m^{\rm Leb}(A)\cdot\frac{e_m}{2^{2n+2}}<m^{\rm Leb}(C_n\cap U)<(\ell^*-
\ell)\cdot\Big(\frac{2^{3n+3}}{e_{n-1}}\Big)^{\ell^*-\ell-1}\cdot e_n,\]
and hence 
\[m^{\rm Leb}(A)<\Big(\frac{2^{3n+3}}{e_{n-1}}\Big)^{\ell^*-\ell-1}\cdot
\frac{2^{2n+2}}{e_{n-1}}\cdot (\ell^*-\ell)\cdot e_n<\Big(\frac{2^{3n+3}}{
e_{n-1}}\Big)^{\ell^*-\ell}\cdot e_n.\]
Pick an open set $U_\ell\subseteq(\can)^{\textstyle\zeta_m}$ such that
$A\subseteq U_\ell$ and $m^{\rm Leb}(U_\ell)<\Big(\frac{2^{3n+3}}{e_{n-1}}
\Big)^{\ell^*-\ell}\cdot e_n$. 

Finally we let $C_{n+1}=C_n\setminus\bigcup\limits_{\ell=1}^{\ell^*}
\big(U_\ell\times(\can)^{\textstyle [\xi_\ell,\zeta)}\big)$. It is easy to
check that $C_{n+1}$ is as required. 
\medskip

After the sets $C_n$ are all constructed we put $C^*=\bigcap\limits_{n<
\omega} C_n$. It follows from $(\circledast)$ that the demand
$(\otimes)^\xi_{C^*}$ is satisfied for each $\xi<\zeta$.
\end{proof}

\begin{lemma}
\label{closincond}
Suppose that $p\in\starQ$, and $q\in\starQ$ is a cover for $p$. Let $C
\subseteq [T^p]\subseteq [T^q]$ be a closed set of positive Lebesgue measure
in $[T^q]$ (so $m^{\rm Leb}(\bh_p[C])>0$). Then there is a condition $p^*\in 
\starQ$ stronger than $p$ and such that $[T^{p^*}]\subseteq C$.
\end{lemma}

\begin{proof}
For $t\in K^*$ let $F_t:[0,1]^{\textstyle\pos(t)}\longrightarrow [0,1]$ be
defined by
\[F_t(r_\nu:\nu\in\pos(t))=\frac{\sum\{r_\nu:\nu\in\pos(t)\}}{2^{N_{k_t}
-|g_t|}}.\] 
This defines a function $\bF$ on $K^*$. Plainly, $(K^*,\Sigma^*,\bF)$ is a 
nice measured tree creating pair (we are going to use it to simplify
notation only). 

Let $T\subseteq T^p$ be a tree such that $\max(T)=\emptyset$ and
$C=[T]$. For $\eta\in T$ let $t_\eta\in \Sigma^*(t^p_\eta)$ be such that 
\[\pos(t_\eta)=\suc_T(\eta),\quad \nor[t_\eta]=\nor[t^p_\eta]\quad\mbox{ and
}\quad g_{t_\eta}=g_{t^p_\eta}.\]
Let $q=\langle t_\eta:\eta\in T\rangle$. It should be clear that (as $C$ has
positive Lebesgue measure) $q$ is a condition in $\bQ^\mtree_4(K^*,\Sigma^*,
\bF)$ (note: $\bF$, not $\bF^*$!).  Moreover, possibly shrinking $T$ and
$C$, we may request that 
\begin{itemize}
\item $\nor[t_\eta]>2$ for all $\eta\in T$,
\item $\mu^\bF(q)>1/2$, and $\mu^\bF_q(\eta)\geq 2^{-2^{\lh(\eta)+1}}$ for
each $\eta\in T$ 
\end{itemize}
(remember \ref{nordense}, or actually its proof). Let $k_0=\lh(\mrot(T))$.  

Fix an integer $k>k_0$ for a moment. Let $A=\{\eta\in T:\lh(\eta)=k\}$ (so
it is a front of $T$). For each $\eta\in T[q,A]$, by downward induction, we
define $s_\eta\in\Sigma^*(t_\eta)$ and a real $a_\eta\in [0,1]$ such that  
\begin{enumerate}
\item[$(\bigstar)_\eta$] \quad $a_\eta\geq \prod\limits_{\ell=\lh(\eta)}^{
k-1}\Big(1-2^{-2^{\ell+3}}\Big)\cdot\mu^\bF_q(\eta)$.
\end{enumerate}
If $\eta\in A$, then we let $a_\eta=1$ (and $s_\eta$ is not defined). 

\noindent Suppose that $a_\nu$ has been defined for all $\nu\in\pos(t_\eta)$ 
so that $(\bigstar)_\nu$ holds. Then 
\[\begin{array}{l}
F_{t_\eta}(a_\nu:\nu\in\pos(t_\eta))\geq\prod\limits_{\ell=\lh(\eta)+1}^{k
-1}\Big(1-2^{-2^{\ell+3}}\Big)\cdot F_{t_\eta}(\mu^\bF_q(\nu):\nu\in\pos(
t_\eta))={}\\
\prod\limits_{\ell=\lh(\eta)+1}^{k-1}\Big(1-2^{-2^{\ell+3}}\Big)
\cdot \mu^\bF_q(\eta)\geq 2^{-2^{\lh(\eta)+3}}
  \end{array}\]
(remember our requests on $q$). Consequently we may apply \ref{Sah7.2} (for
$t=t_\eta$, $r_\nu=a_\nu$ and $g'=g_{t_\eta}$) to pick $s_\eta\in
\Sigma^*(t_\eta)$ such that 
\begin{enumerate}
\item[$(\alpha)$] $\nor[s_\eta]=\nor[t_\eta]-1$, and 
\item[$(\beta)$]  $a_\eta\stackrel{\rm def}{=} F^*_{s_\eta}(a_\nu:\nu\in
\pos(s_\eta))\geq \Big(1-2^{-2^{\lh(\eta)+3}}\Big)\cdot F_{t_\eta}(a_\nu:\nu
\in\pos(t_\eta))\geq \prod\limits_{\ell=\lh(\eta)}^{k-1}\Big(1-2^{-2^{\ell+
3}}\Big)\cdot\mu^\bF_q(\eta)$.
\end{enumerate}
This completes the choice of $s_\eta$'s and $a_\eta$'s. Now we build a
system $\langle s^k_\eta:\eta\in S_k\setminus\max(S_k)\rangle$ such that
$S_k\subseteq T[q,A]$ is a finite tree, $\mrot(S_k)=\mrot(T)$,
$s^k_\eta=s_\eta$ and $\suc_{S_k}(\eta)=\pos(s^k_\eta)$ for $\eta\in S_k
\setminus\max(S_k)$. 

Next, applying K\"onig Lemma, we pick an infinite set $I\subseteq\omega$ and
a system $p^*=\langle t^{p^*}_\eta:\eta\in T^{p^*}\rangle\in\bQ^*_\emptyset(
K^*,\Sigma^*)$ such that $\mrot(T^{p^*})=\mrot(T)$ and 
\[\eta\in T^{p^*}\ \&\ k_1,k_2\in I\ \&\ \lh(\eta)<k_1<k_2\quad \Rightarrow
\quad t^{p^*}_\eta=s^{k_2}_\eta.\]
It follows from our construction that necessarily $p^*\in \starQ$, and it is
a condition stronger than $p$, and $[T^{p^*}]\subseteq [T]=C$.
\end{proof}

\begin{theorem}
In $\bV^{\bP_{\omega_2}}$, the condition (b) of \ref{reduction} holds.
\end{theorem}

\begin{proof}
For $\alpha<\omega_2$ let $\dot{x}_\alpha$ be a $\bP_\alpha$--name for the
generic real added at stage $\alpha$ (so it is a member of $\can$ if
$\alpha\in Z$, and a member of $\prod\limits_{k<\omega}\bH^*(k)$ if $\alpha
\in\omega_2\setminus Z$). 

Suppose that $\dot{f}^*$ is a $\bP_{\omega_2}$--name for a function from
$\can$ to $\can$, and $p\in\bP_{\omega_2}$. 

For each $\delta\in Z$ pick a template with a name $(\bar{\bt}^\delta,
\bar{\tau}^\delta)$, an enumeration $\bar{\zeta}^\delta=\langle
\zeta^\delta_n:n<\omega\rangle$ of $\zeta_{\bar{\bt}^\delta}$, and a
condition $p^\delta\in\bP_{\omega_2}$ such that 
\begin{itemize}
\item $\zeta_{\bar{\bt}^\delta}\geq\omega$, $\bar{\bt}^\delta$ behaves well
for $\bar{\zeta}^\delta$ (see \ref{easytemp}(1)), 
\item $p^\delta\geq p$ and $(p^\delta,\dot{f}^*(\dot{x}_\delta))$ obeys
$(\bar{\bt}^\delta,\bar{\tau}^\delta)$,
\item $\delta\in w^{\bar{\bt}^\delta}$ and $w^{\bar{\bt}^\delta}\setminus
(\delta+1)\neq\emptyset$.  
\end{itemize}
Using Fodor Lemma (and \ref{easycount}(2)) we find a template with a name
$(\bar{\bt},\bar{\tau})$, ordinals $\zeta^*<\zeta_{\bar{\bt}}$ and
$\xi<\omega_2$, an enumeration $\bar{\zeta}=\langle\zeta_n:n<\omega\rangle$
of $\zeta_{\bar{\bt}}$, and a stationary set $Z^* \subseteq Z$ such that for
each $\delta,\delta'\in Z^*$ we have 
\begin{enumerate}
\item[(i)]   $(\bar{\bt}^\delta,\bar{\tau}^\delta)$ is isomorphic to
$(\bar{\bt},\bar{\tau})$ by an isomorphism mapping $\bar{\zeta}^\delta$ to
$\bar{\zeta}$, and $w^{\bt}=\zeta_{\bar{\bt}}$, $\bar{\bt}=\langle 
\bt_n:n<\omega\rangle$, $\bar{\tau}=\langle\tau_n:n<\omega\rangle$, and
\item[(ii)]  $\otp(w^{\bar{\bt}^\delta}\cap\delta)=\zeta^*$,
$w^{\bar{\bt}^\delta}\cap\delta\subseteq\xi$, and $p\in\bP_\xi$ and 
\item[(iii)] $\bar{\bt}^\delta\restriction\xi=\bar{\bt}^{\delta'}
\restriction\xi$. 
\end{enumerate}
Let $\dot{A}$ be the $\bP_{\omega_2}$--name for the set $\{\dot{x}_\delta:
\delta\in Z^*\ \&\ p^\delta\in\Gamma_{\bP_{\omega_2}}\}$ and let $\Psi:(
\can)^{\textstyle [\zeta^*+1,\zeta_{\bar{\bt}})}\longrightarrow\can$ be the 
canonical homeomorphism (induced by a bijective mapping from $\omega\times
[\zeta^*+1,\zeta_{\bar{\bt}})$ onto $\omega$). Now, in  $\bV^{\bP_{
\omega_2}}$, we define a mapping $\dot{f}_1:\dot{A}\longrightarrow\can$ by:  
\[\dot{f}_1(\dot{x}_\delta)=\Psi\Big(\psi^{\bar{\bt}^\delta}_{\zeta_{
\bar{\bt}}}(\dot{x}_\alpha:\alpha\in w^{\bar{\bt}^\delta})\restriction
[\zeta^*+1,\zeta_{\bar{\bt}})\Big)\]
($\psi^{\bar{\bt}^\delta}_{\zeta_{\bar{\bt}}}$ is as defined before
\ref{templates}). Let $p^*=p^\delta\restriction\delta$ for some
(equivalently: all) $\delta\in Z^*$. 
\begin{claim}
\label{cl5}
\[p^*\forces_{\bP_{\omega_2}}\mbox{`` the set }\{(x,\dot{f}_1(x)):x\in
\dot{A}\}\mbox{ has positive outer measure ''.}\]
\end{claim}

\begin{proof}[Proof of the claim]
Assume not. Then there are an ordinal $\xi^*$, a condition $q$, and a
$\bP_{\omega_2}$--name $\dot{D}$ such that  
\begin{itemize}
\item $\xi\leq\xi^*<\omega_2$, $q\in\bP_{\xi^*}$, and $q\geq p^*$, 
\item $\dot{D}$ is a $\bP_{\xi^*}$--name for a (Lebesgue) null subset of
$(\can)^{\textstyle[\zeta^*,\zeta_{\bar{\bt}})}$, and
\item $q\forces_{\bP_{\omega_2}}$`` $(\forall\delta\in Z^*)\big(p^\delta\in
\Gamma_{\bP_{\omega_2}}\ \Rightarrow\ \psi^{\bar{\bt}^\delta}_{\zeta_{
\bar{\bt}}}(\dot{x}_\alpha:\alpha\in w^{\bar{\bt}^\delta})
\restriction [\zeta^*,\zeta_{\bar{\bt}})\in\dot{D})$ ''.  
\end{itemize}
(Note that above we use the fact that the forcing used at $\delta\in Z$ is
the random real forcing, so the conditions are closed sets of positive
measure. This allows us to replace $\big(\dot{x}_\delta,\psi^{\bar{
\bt}^\delta}_{\zeta_{\bar{\bt}}}(\dot{x}_\alpha:\alpha\in w^{
\bar{\bt}^\delta})\restriction [\zeta^*+1,\zeta_{\bar{\bt}})
\big)$ by $\psi^{\bar{\bt}^\delta}_{\zeta_{\bar{\bt}}}(\dot{x}_\alpha:
\alpha\in w^{\bar{\bt}^\delta})\restriction [\zeta^*,\zeta_{\bar{\bt}})$.)
Fix any $\delta^*\in Z^*$ larger than $\xi^*$ and let $\langle\alpha_\zeta: 
\zeta<\zeta^*\rangle$ be the increasing enumeration of $w^{\bar{\bt}^{
\delta^*}}\cap\delta^*$ and let $\dot{z}_\zeta=\dot{x}_{\alpha_\zeta}$,
and $\dot{\bar{z}}=\langle\dot{z}_\zeta:\zeta<\zeta^*\rangle$. Note that the
conditions $p^{\delta^*}$ and $q$ are compatible. Also, as
$\dot{x}_{\delta^*}$ is (a name for) a random real over
$\bV^{\bP_{\delta^*}}$, we have
\[\begin{array}{ll}
q\forces_{\bP_{\delta^*+1}}&\mbox{`` the set}\\
&\quad\dot{B}\stackrel{\rm def}{=}\{\bar{y}\in (\can)^{\textstyle
[\zeta^*+1,\zeta_{\bar{\bt}})}:\big\langle\psi^{\bar{\bt}^{\delta^*}}_{
\zeta^*+1}(\dot{\bar{z}}\conc\langle\dot{x}_{\delta^*}\rangle)(\zeta^*)\big
\rangle\conc\bar{y}\in\dot{D}\}\\
&\ \mbox{ is null ''.}
\end{array} \] 
Using Lemma \ref{specclos}, we may pick (a $\bP_{\delta^*+1}$--name for) a
closed set $\dot{C}^*\subseteq (\can)^{\textstyle [\zeta^*+1,\zeta_{\bar{
\bt}})}$ such that the condition $q$ forces (in $\bP_{\delta^*+1}$): 
\begin{itemize}
\item $\dot{C}^*\subseteq \{\psi^{\bar{\bt}^{\delta^*}}_{\zeta_{\bar{\bt}}}(
\bar{z})\restriction [\zeta^*+1,\zeta_{\bar{\bt}}):\dot{\bar{z}}\conc\langle
\dot{x}_{\delta^*}\rangle\vartriangleleft\bar{z}\in\cZ^{\bar{\bt}^{
\delta^*}}_{\zeta_{\bar{\bt}}}\}$,
\item $\dot{C}^*\cap \dot{B}=\emptyset$, and 
\item the condition $(\otimes)^\xi_{\dot{C}^*}$ of \ref{specclos} holds true
for every $\xi\in [\zeta^*+1,\zeta_{\bar{\bt}})$.
\end{itemize}
(For the first demand remember that $\bar{\bt}^{\delta^*}$ is well behaving,
so the set on the right-hand side has positive Lebesgue measure.) But now,
using \ref{closincond}, we may inductively build a condition $q'\in
\bP_{\omega_2}$ stronger than both $q$ and $p^{\delta^*}$ (and with
the support included in $(\delta^*+1)\cup w^{\bar{\bt}^{\delta^*}}$) such
that  
\[q'\forces_{\bP_{\omega_2}}\mbox{`` }\psi^{\bar{\bt}^{\delta^*}}_{
\zeta_{\bar{\bt}}}(\dot{x}_\alpha:\alpha\in w^{\bar{\bt}^{\delta^*}})
\restriction [\zeta^*+1,\zeta_{\bar{\bt}})\notin\dot{B}\mbox{ '',}\]  
getting an immediate contradiction. 
\end{proof}

Pick any $\delta^*\in Z^*$ and let $\dot{\bar{z}}=\langle\dot{z}_\zeta:
\zeta<\zeta^*\rangle$ be as defined in the proof of \ref{cl5} above. Let
$\dot{E}$ be a $\bP_{\delta^*}$--name for the set
\[\{(r_0,r_1)\in\can\times\can:\dot{\bar{z}}\conc\langle r_0\rangle\in
\cZ^{\bar{\bt}}_{\zeta^*+1}\mbox{ and }\psi^{\bar{\bt}}_{\zeta^*+1}(
\dot{\bar{z}}\conc\langle r_0\rangle)\conc\Psi^{-1}(r_1)\in\rng(
\psi^{\bar{\bt}}_{\zeta^{\bar{\bt}}})\}.\]
So $\dot{E}$ is (a name for) a closed subset of $\can\times\can$. Let
$\dot{f}_2$ be a name for a Borel function from $\can\times\can$ to $\can$
such that 
\begin{quotation}
\noindent if $(r_0,r_1)\in\dot{E}$, and $\psi^{\bar{\bt}}_{\zeta^*+1}(
\dot{\bar{z}}\conc\langle r_0\rangle)\conc\Psi^{-1}(r_1)=\psi^{\bar{
\bt}}_{\zeta^{\bar{\bt}}}(\langle z_\zeta:\zeta<\zeta^{\bar{\bt}}\rangle)$, 

\noindent then for each $n<\omega$
\[\dot{f}_2(r_0,r_1)\restriction n=\tau_n(z_\zeta\restriction \bk^{\bt_n}(
\zeta):\zeta\in w^{\bt_n})\]
\end{quotation}
(remember (i)). It should be clear that $\dot{f}_2$ is (a name for) a
continuous function and  
\[p^*\forces_{\bP_{\omega_2}}\mbox{`` }(\forall x\in\dot{A})(\dot{f}^*(x)=
\dot{f}_2(x,\dot{f}_1(x)))\mbox{ '',}\]
finishing the proof. 
\end{proof}

\begin{corollary}
\label{twomain}
It is consistent that 
\begin{itemize}
\item every sup-measurable function is Lebesgue measurable, and 
\item for every function $f:\mbR\longrightarrow\mbR$ there is a continuous 
function $g:\mbR\longrightarrow\mbR$ such that the set $\{x\in\mbR:f(x)
=g(x)\}$ has positive outer measure.
\end{itemize}
\end{corollary}


\end{document}